\documentclass[10pt]{article}
\usepackage{amsbsy,amssymb, amsmath}
\usepackage{caption,graphicx}
\usepackage[text={6.75in,9.5in},centering,letterpaper]{geometry}
\usepackage[numbers,comma,square,sort&compress]{natbib}
\usepackage{pstricks}

\setcaptionmargin{0.25in}

\newtheorem{Lemma}{Lemma}[section]

\newtheorem{Hypothesis}[Lemma]{Hypothesis}
\newtheorem{Proposition}[Lemma]{Proposition}
\newtheorem{Remark}[Lemma]{Remark}
\newtheorem{Theorem}{Theorem}

 {\begin{trivlist}\item[]\textbf{Proof#1 }}%
 {\hspace*{\fill}$\rule{0.3\baselineskip}{0.35\baselineskip}$\end{trivlist}}

\makeatletter\@addtoreset{equation}{section}\makeatother

\def\bexe{\begin{exercise}}\def\eexe{\eex\end{exercise}}
\def\bsol{\begin{solution}}\def\esol{\eex\end{solution}}
\def\bexa{\begin{example}}\def\eexa{\eex\end{example}}
\def\brem{\begin{remark}}\def\erem{\eex\end{remark}}

\newcommand{\R}{{\mathbb R}}\newcommand{\C}{{\mathbb C}}
\newcommand{\N}{{\mathbb N}}
\newcommand{\Z}{{\mathbb Z}}

\def\CT{{\cal T}}
  
\def\CF{{\cal F}}\def\CG{{\cal G}}\def\CL{{\cal L}}
\def\CO{{\cal O}}\def\CR{{\cal R}}\def\CS{{\cal S}}
\def\CM{{\cal M}}\def\CN{{\cal N}}\def\CB{{\cal B}}

\def\ga{\gamma}\def\om{\omega}\def\th{\theta}\def\uh{\hat{u}}

\def\vt{\vartheta}\def\pa{{\partial}}\def\lam{\lambda}

\newcommand{\bi}{\begin{itemize}}\newcommand{\ei}{\end{itemize}}
\newcommand{\bce}{\begin{center}}\newcommand{\ece}{\end{center}}
\newcommand{\reff}[1]{(\ref{#1})}

\newcommand{\vs}[1]{{\vspace{#1}}}
\def\eps{\varepsilon}
\def\ra{\rightarrow}

\newcommand{\barr}{\begin{array}}\newcommand{\earr}{\end{array}}

\newcommand{\bpm}{\begin{pmatrix}}\newcommand{\epm}{\end{pmatrix}}
\newcommand{\ba}{\begin{array}}\newcommand{\ea}{\end{array}}
\def\dd{\, {\rm d}}\def\ri{{\rm i}}
\def\res{{\rm Res}}\def\er{{\rm e}}
\def\re{{\rm Re}}\def\im{{\rm Im}}
\def\om{\omega}
\def\omnl{\omega_\mathrm{nl}}
\def\hot{{\rm h.o.t}}
\def\del{\delta}

\def\eex{\hfill\mbox{$\rfloor$}}

\def\sig{\sigma}
\def\al{\alpha}
\def\fti{\tilde{f}}\def\gti{\tilde{g}}\def\wti{\tilde{w}}

\def\kap{\kappa}

\def\bd{\begin{displaymath}} \def\ed{\end{displaymath}}
\def\ba{\begin{array}} \def\ea{\end{array}}

\def\eps{\varepsilon}

\def\qed{\hfill $\Box$}\def\rmi{{\rm i}}

\def\uti{\tilde{u}}\def\vti{\tilde{v}}
\def\hti{\tilde{h}}

\def\erf{{\rm erf}}
\def\uhat{\hat{u}}\def\qhat{\hat{q}}
\def\CJ{{\cal J}}\def\CX{{\cal X}}

\def\omnl{\om}
\def\rmd{\dd}\def\rme{\er}\def\rmO{\CO}
\def\uic{u_{{\rm ic}}}

\def\CLt{\tilde{\CL}}

\def\ati{\tilde{\al}}\def\Lamti{\tilde{\Lambda}}\def\lamti{\tilde{\lam}}
\def\uad{u_\mathrm{ad}}\def\CNti{\tilde{\CN}}

\def\philim{\phi_{{\rm lim}}}\def\psilim{\psi_{{\rm lim}}}
\def\psilimti{\tilde{\psi}_{{\rm lim}}}
\def\supp{{\rm supp}}

\def\pmfsti{\tilde{P}_{{\rm mf}}^{{\rm s}}}
\def\etati{\tilde{\eta}}
\def\Kti{{\tilde K}}
\def\Hhat{\hat{H}}\def\Uti{\tilde{U}}\def\Gti{\tilde{G}}
\def\Vti{\tilde{V}}\def\Wti{\tilde{W}}
\def\alti{\tilde{\al}}
\def\ds{\displaystyle}

\def\Re{\mathop\mathrm{Re}\nolimits}    % real part
\def\huga#1{\begin{gather} #1 \end{gather}}

\begin{document}

\title{Diffusive mixing of periodic wave trains in reaction-diffusion systems}

\author{%
Bj\"orn Sandstede\\
Division of Applied Mathematics\\
Brown University\\
Providence, RI~02912, USA
\and
Arnd Scheel\\
School of Mathematics\\
University of Minnesota\\
Minneapolis,  MN 55455, USA
\and
Guido Schneider\\
Institut f\"ur Analysis, Dynamik und Modellierung\\
Universit\"at Stuttgart\\
D-70569 Stuttgart, Germany
\and
Hannes Uecker\\
Institut f\"{u}r Mathematik\\
Universit\"{a}t Oldenburg\\
D-26111 Oldenburg, Germany
}

\date{\today}
\maketitle

\begin{abstract}
We consider reaction-diffusion systems on the infinite line that exhibit a family of spectrally stable spatially periodic wave trains $u_0(kx-\om t;k)$ that are parameterized by the wave number $k$. We prove stable diffusive mixing of the asymptotic states $u_0(k x+\phi_{\pm};k)$ as $x\ra \pm\infty$ with different phases $\phi_-\neq\phi_+$ at infinity for solutions that initially converge to these states as $x\ra \pm\infty$. The proof is based on Bloch wave analysis, renormalization theory, and a rigorous decomposition of the perturbations of these wave solutions into a phase mode, which shows diffusive behavior, and an exponentially damped remainder. Depending on the dispersion relation, the asymptotic states mix linearly with a Gaussian profile at lowest order or with a nonsymmetric non-Gaussian profile given by Burgers equation, which is the amplitude equation of the diffusive modes in the case of a nontrivial dispersion relation.
\end{abstract}

\newpage

%%%%%%%%%%%%%%%%%%%%%%%%%%%%%%%%%%%%%%%%%%%%%%%%%%%%%%%%%%%%%%%%%%%%%%%%%%%%

\section{Introduction}\label{i-sec}

We consider spatially extended pattern-forming systems that exhibit
periodic travelling-wave solutions $u(x,t)=u_0(kx-\om t;k)$ for a
certain range of wave numbers $k\in(k_l,k_r)$.  The profile
$u_0(\th;k)$ is assumed to be $2\pi$-periodic in $\th=kx-\om t$, where
the wave number $k$ and the temporal frequency $\om$ are assumed to be
related via a nonlinear dispersion relation $\om=\omnl(k)$. Examples
are the Taylor vortices in the Taylor--Couette problem, roll solutions
in convection problems, or periodic wave trains in reaction-diffusion
systems.

\begin{figure}
\centering
%\vs{100mm}
%\includegraphics[height=100mm]{./pdf/figure-1}
\includegraphics[height=100mm]{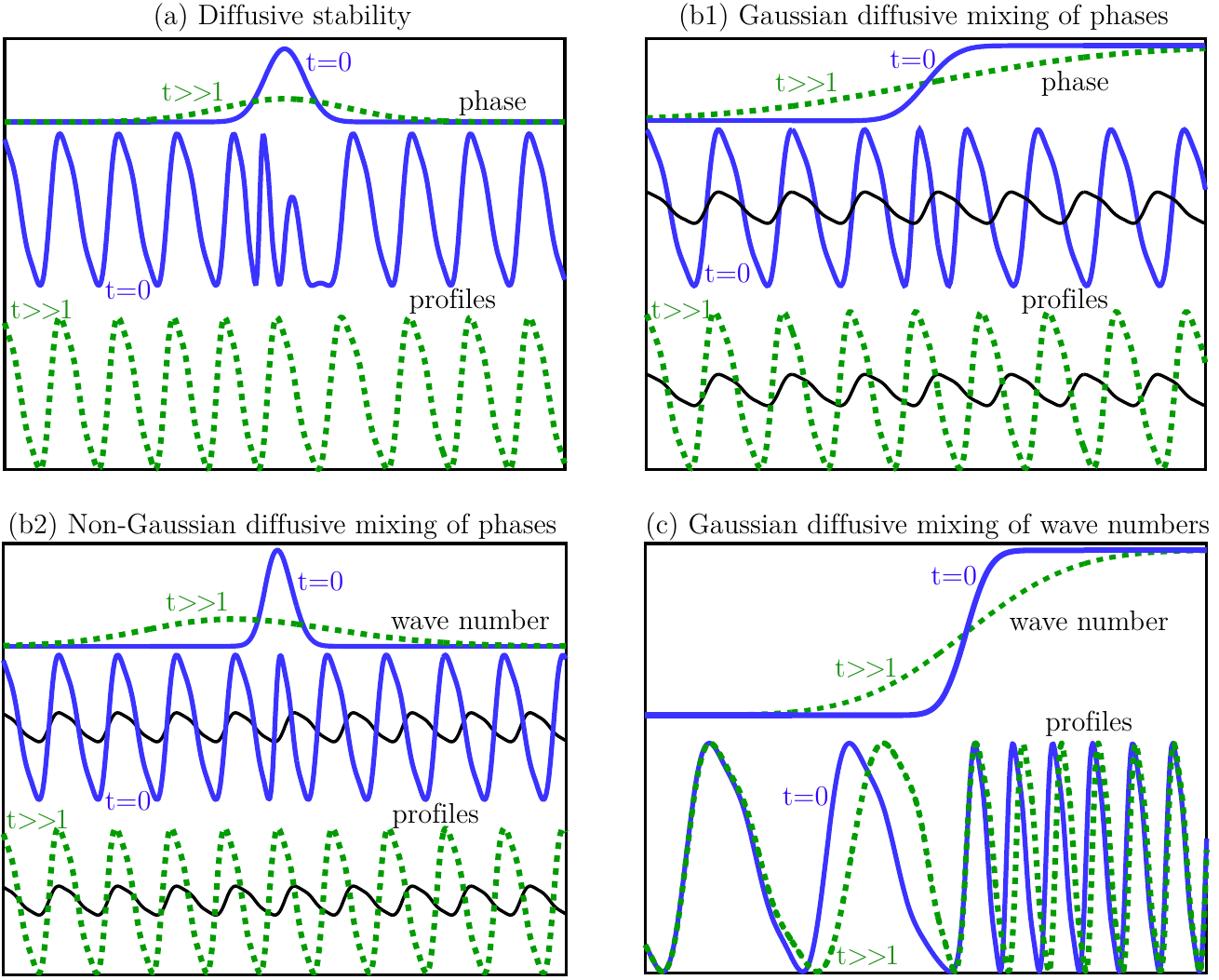}
\caption{The panels illustrate different types of diffusive behaviour in the frame that moves with the speed of the group velocity $c_g$.
(a) Localized perturbations of wave trains decay diffusively like Gaussians;
(b1) If $\omega''=0$, then phase fronts develop when $\phi_-\neq\phi_+$, and the wave number perturbation decays diffusively like a Gaussian;
(b2) If $\omega''\neq 0$, then phase fronts develop, and the wave number perturbation decays as determined by the Burgers equation;
(c) Shown is the expected formal diffusive mixing of wave-number fronts in case $\omega''=0$.
Solid and dashed lines indicate solutions at $t=0$ and for $t\gg1$, respectively, while the small solid curves in (b1)-(b2) indicate the amplitude-scaled spatially periodic wave train to visualize the phase shifts.
\label{blfig1}}
\end{figure}

We are interested in the dynamics of perturbations of wave-train solutions of the above form. Since the linearization around a wave train always possesses essential spectrum up to the imaginary axis, we cannot expect exponential relaxation towards the original profile even for spectrally stable wave trains. Moreover, the periodic nature of the underlying wave train suggests that we should allow perturbations that change the phase or the wave number of the underlying profile. In these cases, we expect that diffusive decay or diffusive mixing of phases or wave numbers dominate the dynamics. In more detail, given a spatially periodic wave train $u_0(k_0x-\om_0t;k_0)$ we may consider (a) its diffusive stability, that is, its stability with respect to spatially localized perturbations, or else the diffusive mixing of the asymptotic states $u_0(k_{\pm}x+\phi_{\pm};k_{\pm})$ as $x\ra \pm\infty$ with (b) identical wave number $k_-=k_+$ but different phases $\phi_-\neq\phi_+$ or (c) 
 different wave numbers $k_-\neq k_+$ for initial data that converge to these states as $x\ra \pm\infty$. More precisely, consider an initial condition of the form
\begin{equation}\label{e:ic}
u(x,0) = u_0(q_0(x)x+\phi_0(x);q_0(x)), \qquad q_0(x)\to k_\pm,\; \phi_0(x)\to\phi_\pm \mbox{ as } x\to\pm\infty,
\end{equation}
where the functions $q_0(x)$ and $\phi_0(x)$ are bounded and small in an appropriate norm. We may then expect that the solution $u(t,x)$ can, to leading order, be written in the form
\[
u(x,t) \approx u_0(q(x,t)x+\phi(x,t)-\omega_0 t;q(x,t)),
\]
and the issue is to determine the behaviour of the phase $\phi(x,t)$ and the local wave number $q(x,t)$ as $t\to\infty$. As indicated above, we can  distinguish three different classes of initial data, namely (a) constant wave number $q_0(x)\equiv k_0$ and equal phases $\phi_+=\phi_-$ at infinity for non-zero phase perturbations $\phi_0(x)\not\equiv0$, which correspond to localized perturbations of the underlying wave train, (b) constant wave number $q_0(x)\equiv k_0$ but different phases $\phi_+\neq\phi_-$ at infinity, which correspond to a relative phase shift of the wave train at $\pm\infty$, and (c) different wave numbers $k_-\neq k_+$ at infinity; see Figure~\ref{blfig1} for an illustration.

In this paper, we address the cases (a) and (b) for general reaction-diffusion systems
\begin{gather}
\pa_t u=D\pa_x^2 u+f(u)
\label{rd1}
\end{gather}
with $x\in\R$, $t\geq0$, and $u(x,t)\in\R^d$, where $D\in\R^{d\times d}$ is symmetric and positive definite, and $f$ is smooth. We now outline our results and refer to Theorems~\ref{thstabil} and~\ref{thmix1} for the precise statements:
\bi
\item[(a)]
For \textbf{localized perturbations} of a single wave train, that is, for $q_0(x)\equiv k_0$ and $\phi_-=\phi_+$, we transfer existing stability results from specific systems \cite{gs96a,gs98,gstohoku,ue04nsb} to general reaction-diffusion systems. In lowest order, the dynamics near a wave train can be described by the evolution of the local wave number $q(t,x)$, and we prove that the renormalized wave number difference $t[q(t^{1/2}x,t)-k_0]$ converges towards a multiple of 
the $x$-derivative of the Gaussian $\frac{1}{\sqrt{4\pi\al}}\exp(-\frac{x^2}{4\al})$ for an appropriate constant $\alpha>0$. This yields the asymptotics
\[
\sup_{x\in \R}|u(x,t)-\phi_{\rm lim}\phi^*(x-c_g t,t) \pa_\th u_0(\th;k)| \leq C_2 t^{-1+b}\text{ as } t\ra\infty, 
\]
with $\phi^*(x,t)=\frac{1}{\sqrt{4\pi\al t}}\mathrm{e}^{-x^2/(4 \al t)}$, where $\phi_{\rm lim}\in\R$ depends on the initial data, where $\al>0$ and $c_g\in\R$ are constants determined by the spectral properties of $u_0(\cdot,k_0)$, and where $b>0$ is a small, but arbitrary, correction coefficient.
\item[(b)]
For perturbation that induce a \textbf{global phase shift}, that is, for $q_0(x)\equiv k_0$ and $\phi_-\neq\phi_+$ but 
with $|\phi_d|:=|\phi_+-\phi_-|$ small, we establish diffusive decay of wave-number perturbations. Specifically, the renormalized wave number converges to a Gaussian profile when $\om''(k_0)=0$, while it converges to a nonsymmetric non-Gaussian profile when $\om''(k_0)\neq0$. This latter case is the major result of this paper.
\ei

The case (c) where $q_0(x)\to k_\pm$ as $x\to\pm\infty$ with $k_-\neq k_+ $ is more difficult and depends crucially on the sign of $\omega''(k_0)$. If $\om''(k_0)\neq0$, diffusive mixing of the local wave number cannot be expected: instead, depending on the sign of $\om''(k_0)(k_+-k_-)$, we expect that $q(x,t)$ evolves either as a stable viscous shock or as an approximate rarefaction wave \cite{DSSS05}. If $\omega''(k_0)=0$, nonlinear diffusive mixing can be expected, but, for some technical issues that we explain below, a rigorous proof remains open and is left for future research. 

The proof of diffusive mixing of phases of wave trains in systems with no $S^1$-symmetry has resisted many attempts. With the rigorous separation of the phase variable $\phi$ from remaining modes found in \cite{DSSS05}, a new technique is now available to treat this question. This method combined with the renormalization group method \cite{bk92,bkl94}, which has been applied for instance in \cite{gs96a,gs98,es00b,gsu04,ue04ibl,ue04nsb} to a variety of pattern-forming and hydrodynamic systems, finally yields our results. Diffusive mixing results for the real Ginzburg--Landau equation, which has a natural decomposition into phase and amplitude variables due to its gauge symmetry, have been obtained for instance in \cite{bk92,gm98}.

The results in this paper were presented at the Snowbird meeting in 2007. 
Meanwhile, similar results on the diffusive stability of wave trains have 
been established in \cite{jz10,jnrz11} using pointwise estimates.

\paragraph{Notation.}

Throughout this paper, we denote many different constants that are
independent of the Burgers parameters $\al,\beta$ and the rescaling
parameter $L>0$ by the same symbol $C$. For $m_1,m_2\in\N$, we define
the weighted spaces $H^{m_2}(m_1){=}\{u\in L^2(\R) : \| u
\|_{H^{m_2}(m_1)}{<}\infty\}$ with norm $\| u\|_{H^{m_2}(m_1)} = \| u
\rho^{m_1} \|_{H^{m_2}(\R)}$, where $ \rho(x) = (1+x^2)^{1/2}$ and
$H^{m_2}(\R)$ is the Sobolev space of functions with weak derivatives
up to order $m_2$ in $L^2(\R)$. With an abuse of notation, we
sometimes write $\|u(x,t)\|_{H^{m_2}(m_1)}$ for the
$H^{m_2}(m_1)$-norm of the function $x\mapsto u(x,t)$. The Fourier
transform is denoted by $\CF$ so that $\uh(k):=\CF(u)(k)=
\frac{1}{2\pi} \int \er^{-\ri k x} u(x) \dd x$ for $u\in
L^2(\R)$. Parseval's identity and $\CF(\pa_x u)(k)=\ri k\uh(k)$ imply
that $\CF$ is an isomorphism between $H^{m_2}(m_1)$ and
$H^{m_1}(m_2)$, that is, the weight in physical space yields
smoothness in Fourier space and vice versa. To indicate functions in
Fourier space, we also write $\uhat\in\Hhat^{m_1}(m_2)$ instead of
$\uhat\in H^{m_1}(m_2)$.

%%%%%%%%%%%%%%%%%%%%%%%%%%%%%%%%%%%%%%%%%%%%%%%%%%%%%%%%%%%%%%%%%%%%%%%%%%%%

\section{Statement of results}

\subsection{Wave trains and their dispersion relations}\label{is-sec}

We assume that there are numbers $k_0\neq0$ and $\omega_0\in\R$ such that \reff{rd1} has a  solution of the form $u(x,t)=u_0(k_0x-\om_0 t)$, where $u_0(\th)$ is $2\pi$-periodic in its argument. Thus, $u_0$ is a $2\pi$-periodic solution of the boundary-value problem
\begin{equation}\label{bvpnl}
k^2 D \pa_{\th}^2 u+ \om \pa_\th u + f(u) = 0 
\end{equation}
with $k=k_0$ and $\om=\om_0$. Linearizing \reff{bvpnl} at $u_0$ yields the linear operator 
\begin{equation}\label{bvpl}
\CL_0=\CL(k_0)=k_0^2 D \pa_{\th}^2 + \om_0 \pa_\th + f^\prime(u_0(\th)),
\end{equation}
which is closed and densely defined on $L^2_\mathrm{per}(0,2\pi)$ with domain $\mathcal{D}(\CL_0)=H^2_\mathrm{per}(0,2\pi)$. We assume that $\lambda=0$ is a simple eigenvalue of $\CL_0$ on $L^2_\mathrm{per}(0,2\pi)$, so that its null space is one-dimensional and therefore spanned by the derivative $\pa_\th u_0$ of the wave train.

We may now vary the parameter $k$ in (\ref{bvpnl}) near $k=k_0$ and again seek $2\pi$-periodic solutions of (\ref{bvpnl}). The derivative of the boundary-value problem (\ref{bvpnl}) with respect to $\om$, evaluated at $k=k_0$ in the solution $u_0$, is given by $\pa_\th u_0$. Since $\lambda=0$ is a simple eigenvalue of $\CL_0$ on $L^2_\mathrm{per}(0,2\pi)$, we see that $\pa_\th u_0$ does not lie in the range of $\CL_0$, and the linearization of the boundary-value problem (\ref{bvpnl}) with respect to $(u,\om)$ is therefore onto. Thus, exploiting the translation symmetry of (\ref{bvpnl}) we can solve (\ref{bvpnl}) uniquely, up to translations in $\th$, for $(u,\om)$ as functions of $k$ and obtain the wave trains
\begin{equation}\label{nld}
u(x,t) = u_0(kx-\omnl(k)t;k), \quad k\in(k_l,k_r),
\end{equation}
where $\omnl(k_0)=\om_0$ and $k_l<k_0<k_r$. In particular, wave trains exist for wave numbers $k$ in an open interval centered around $k_0$. We call the function $k\mapsto\omnl(k)$ the nonlinear dispersion relation and define the phase speed of the wave train with wave number $k$ by $c_\mathrm{p}:=\omnl(k)/k$ and its group velocity by
\begin{equation}\label{gv}
c_\mathrm{g} = \frac{\rmd\omnl}{\rmd k}(k).
\end{equation}
To state our assumptions on the spectral stability of the wave train $u_0$ as a solution to the reaction-diffusion system (\ref{rd1}), we consider the linearization 
\begin{equation}\label{rdlc}
\pa_t v = \CL_0 v
\end{equation}
of \reff{rd1} in the frame $\th=k_0x-\om_0 t$ that moves with the phase speed $c_\mathrm{p}=\om_0/k_0$. Particular solutions to this problem can be found through the Bloch-wave ansatz
\begin{equation}\label{bwa}
u(\th,t) = \rme^{\lam(\ell)t+\ri\ell\th/k_0} \vti(\th,\ell), 
\end{equation}
where $\ell\in\R$ and $\vti(\th,\ell)$ is $2\pi$-periodic in $\th$ for each $\ell$. In fact, since $\tilde{v}(\vt,\ell+k_0) = \rme^{\rmi\vt} \tilde{v}(\vt,\ell)$, we can restrict $\ell$ to the interval $[-k_0/2,k_0/2)$. Substituting \reff{bwa} into \reff{rdlc}, we obtain 
\begin{gather}
\CLt(\ell)\vti=\lam(\ell)\vti
\label{levp}
\end{gather}
with a family of operators $\CLt$ given by
\begin{align}
\CLt(\ell)\vti &= \er^{-\ri \ell\th/k_0}\CL(\er^{\ri\ell\th/k_0}\vti(\th,\ell))
= k_0^2 D \left(\pa_\th+\ri\ell/k_0\right)^2 \vti + \om \left(\pa_\th+\ri\ell/k_0\right) \vti + f^\prime(u_0(\th)) \vti,
\label{fam}
\end{align}
each of which is a closed operator on $L^2_\mathrm{per}(0,2\pi)$ with dense domain $H^2_\mathrm{per}(0,2\pi)$. In particular, $\CLt(\ell)$ has compact resolvent, and its spectrum is therefore discrete. We can label the eigenvalues of $\CLt(\ell)$ by indices $j\in\N$ and write them as continuous functions $\lambda_j(\ell)$ of $\ell$. In addition, we can order these eigenvalues so that $\Re\lambda_{j+1}(0)\leq\Re\lambda_j(0)$ for all $j$. In fact, the curves $\ell\mapsto\lambda_j(\ell)$ are analytic except possibly near a discrete set of values of $\ell$ where the values of two or more curves $\lambda_j(\ell)$ for different indices $j$ coincide.

Next, we assume that $\lambda_1(0)$ is the rightmost element in the spectrum for $\ell=0$. Since we assumed that $\lambda=0$ is algebraically simple as an eigenvalue of $\CL$, there is a curve $\lambda_1(\ell)$ of eigenvalues with $\lam_1(0)=0$, and this curve is analytic in $\ell$ for $\ell$ close to zero. We call 
the curve $\lam_1(\ell)$ the linear dispersion relation and denote the associated eigenfunctions of $\CLt(\ell)$ by $\vti_1(\th,\ell)$. We shall compute the derivative $\rmd\lambda_1/\rmd\ell$ and recover the group velocity as defined via the nonlinear dispersion relation, namely 
\begin{gather}
-\im\pa_\ell\lam_1|_{\ell=0} = -c_\mathrm{p}+\pa_k\omnl(k_0) = -c_\mathrm{p}+c_\mathrm{g}.
\label{cglin}
\end{gather}
We remark that the phase velocity $c_\mathrm{p}$ appears in this formula solely because we computed $\lam_1$ in the frame moving with speed $c_\mathrm{p}$, while $\omnl$ was computed in the steady frame. We also note that the signs of the second derivatives of $\lam_1$ and $\omnl$ are, in general, not related. Finally, we assume that $\Re\lambda_1^{\prime\prime}(0)<0$ and that all other eigenvalues $\lam_j(\ell)$ satisfy $\Re\lambda_j(\ell)<-\sig_0$. The following hypothesis summarizes the assumptions we made so far.

\begin{Hypothesis}[Existence of spectrally stable wave
  trains]\label{hyp1}
  Equation \reff{rd1} admits a spectrally stable wave train solution
  $u(x,t)=u_0(\th)$ with $\th=k_0x-\om_0 t$ for appropriate numbers
  $k_0\neq0$ and $\omega_0\in\R$, where $u_0$ is
  $2\pi$-periodic. Spectral stability entails the following
  properties. First, the linearization $\CL_0$ of \reff{rd1} about
  $u_0$ has a simple eigenvalue at $\lambda=0$. Furthermore, the
  linear dispersion relation $\lambda_1(\ell)$ with $\lambda_1(0)=0$
  is dissipative so that $\lambda_1^{\prime\prime}(0)<0$, and there
  exist constants $\sigma_0,\ell_0,\al_0>0$ such that
  $\Re\lam_1(\ell)<-\sigma_0$ for $|\ell|>\ell_0$ and
  $\Re\lam_1(\ell)<-\al_0\ell^2$ for $|\ell|<\ell_0$, while all other
  eigenvalues $\lam_j(\ell)$ with $j\ge2$ have $\Re\lam_j(\ell)\le
  -\sigma_0$ for all $\ell\in [-k/2,k/2)$.
\end{Hypothesis}

Standard perturbation theory yields that the wave trains $u_0(kx-\omnl(k)t;k)$ are also spectrally stable, possibly for a smaller interval $k_l\leq \tilde{k}_l<k<\tilde{k}_r\leq k_r$ of wave numbers than the interval of existence. By changing $k_l,k_r$ accordingly, we shall assume from now on that the wave trains $u_0(\cdot;k)$ with $k\in(k_l,k_r)$ are spectrally stable with uniform constants $\ell_0,\sig_0,\al_0$.

For later use, we collect a few properties of the linear dispersion relation and refer to \cite[\S4.2]{DSSS05} for their derivation. We denote by
\[
\CL_\mathrm{ad} u = k_0^2 D \pa_\th^2 u - \om_0 \pa_\th u + f^\prime(u_0(\th))^T u
\]
the $L^2_{{\rm per}}((0,2\pi))$-adjoint of $\CL_0=\CL(k_0)$ and let $u_\mathrm{ad}$ be a nontrivial function in its null space with the normalization
\begin{equation}\label{adno}
\langle u_\mathrm{ad},\pa_\th u_0 \rangle_{L^2(0,2\pi)} = 1.
\end{equation}
Using the adjoint eigenfunction, we have
\begin{align}\label{der}
\lambda_1^\prime(0) & = 
\ri\left\langle u_\mathrm{ad},c_\mathrm{p}\pa_\th u_0
+ 2 k_0 D \pa_{\th}^2 u_0 \right\rangle_{L^2} 
= \ri (c_\mathrm{p}- c_\mathrm{g})\in\ri\R,\\
\lambda_1^{\prime\prime}(0) & = 
-\left\langle u_\mathrm{ad},4 k_0 D \pa_{k}\pa_\th u_0
+ 2 D \pa_\th u_0 \right\rangle_{L^2}\in\R.
\end{align}
We shall also use the identity
\begin{gather}
\pa_\ell \vti_1(\cdot,0)=\ri\pa_k u_0
\label{dlv}
\end{gather}
that was established in \cite[\S4.2]{DSSS05}.

%%%%%%%%%%%%%%%%%%%%%%%%%%%%%%%%%%%%%%%%%%%%%%%%%%%%%%%%%%%%%%%%%%%%%%%%%%%%

\subsection{Statement of results}

Throughout this section, we fix the wave number $k_0$ of a wave train $u_0(k_0 x-\omega_0 t;k_0)$ of the reaction-diffusion system \reff{rd1} that satisfies Hypothesis~\ref{hyp1}. We then set
\[
\om_0=\om(k_0), \quad 
c_p=\om_0/k_0, \quad
c_g=\om'(k_0), \quad
\beta=-\frac 1 2\om''(k_0), \quad 
\th=k_0 x-\om_0 t
\]  
and write
\begin{gather}
\lam_1(\ell) = \ri (c_p-c_g)\ell - \al\ell^2 + \CO(\ell^3)
\label{disp-rem}
\end{gather}
for the expansion of the linear dispersion relation of $u_0(\cdot;k_0)$. 
For convenience henceforth we write $k=k_0$. Before we state our result, we remark that the decomposition of the initial data that we shall use below in the statements of our theorems is not unique. This non-uniqueness will be removed in the proofs but does not affect the conclusions made in the results below.

Our first result states that $u_0$ is diffusively stable with respect to localized perturbations and extracts the leading-order behaviour of the displacement for large times. For notational convenience, in the following we consider 
initial conditions at $t=1$.

\begin{Theorem}[Diffusive stability] \label{thstabil}
Let $u_0(\cdot;k)$ be a spectrally stable wave train that satisfies 
Hypothesis~\ref{hyp1} and pick $b\in(0,1/2)$; then there are 
$\eps,C>0$ such that 
the following holds. If, for some $\th_0\in [0,2\pi)$,
\begin{gather}
u(x,t)|_{t=1} = u_0(\th-\th_0+\phi_0(x);k) + v_0(x) \quad\mbox{with}\quad 
\|\phi_0\|_{H^3(3)}, \|v_0\|_{H^2(3)} \le \eps, 
\label{icdecom1}
\end{gather}
then the solution $u(x,t)$ of \reff{rd1} exists for all times 
$t\geq 1$, it can be written as
\[
u(x,t)=u_0(\th-\th_0+\phi(x,t);k)+v(x,t), 
\]
and there is a constant $\philim\in\R$ depending only on 
the initial condition so that
\begin{equation} \label{n27}
\sup_{x\in\R} \bigl|\phi(x,t)-\philim G(x-c_g(t-1),t)\bigr|+|v(x,t)| 
\leq C t^{-1+b},
\end{equation}
where
\begin{gather}
G(x,t)=\frac{1}{\sqrt{4\al\pi t}}\mathrm{e}^{-x^2/(4\al t)}.
\label{n31}
\end{gather}
In particular, we have
\[
\sup_{x\in\R} |u(x,t)- u_0(\th-\theta_0+\philim G(x-c_g(t-1),t);k)| 
\leq C_1 t^{-1+b}.
\]
\end{Theorem}

Next, we discuss diffusive mixing of phases for non-localized phase 
perturbations. In this situation, the precise asymptotics of 
perturbations depends on $\beta=-\frac 1 2\omega''(k)$.

\begin{Theorem}[Diffusive mixing of phases] \label{thmix1}
Let $u_0(\cdot;k)$ be a spectrally stable wave train that satisfies 
Hypothesis~\ref{hyp1} and pick $b\in(0,1/2)$; 
then there are constants $\eps,C>0$ such that 
the following holds. 

(i) Assume that $\beta=-\frac 1 2\omnl''(k)=0$ and 
$u(x,t)|_{t=1}=u_0(\th-\th_0+\phi_0(x);k)+v_0(x)$ with 
$\phi_0(x)\ra\phi_{\pm}$ for $x\ra\pm\infty$, 
$|\phi_d|=|\phi_+-\phi_-|\le \eps$, and 
\begin{gather}
\|\phi_0'(\cdot)\|_{H^2(2)}, 
\| v_0 \|_{H^2(2)} \le \eps. 
\label{icdecom2}
\end{gather}
Then  the solution $u(x,t)$ to \reff{rd1} exists for all $ t \geq 1 $, 
and can be written as 
$$
u(x,t)=u_0(\th-\th_0+\phi^*(x-c_g(t-1),t);k)+ v(x,t), 
$$
where 
\begin{gather} \label{n31b}
\phi^*(x,t)=\phi_-+(\phi_+-\phi_-)\erf(x/\sqrt{\al t}),\quad \text{ with } 
\quad\erf(x)=\frac 1 {\sqrt{4\pi}} \int_{-\infty}^x \er^{-\xi^2/4}\dd \xi. 
\end{gather}
and $ \sup_{x\in \R} |v(x,t)| \leq C \ t^{-1/2+b}$. 

(ii) 
The same result holds if $\beta=-\frac 1 2\om''(k)\neq 0$, with 
$\phi^*(x,t)$ replaced by 
\begin{gather}
\phi^*(x,t)=\frac{\al}{\beta}
\ln(1+z\erf(x/\sqrt{\al t})), \quad \ln(1+z)=\phi_+-\phi_-. 
\label{n31c}
\end{gather}
\end{Theorem}

\begin{Remark}{\rm 
Clearly, the decompositions \reff{icdecom1} and \reff{icdecom2} are not unique. 
For instance, $\phi_0\equiv 0$ would be one possibility 
in \reff{icdecom1}, 
but we may shift perturbations between $\phi_0$ and $v_0$. 
In the proof we shall fix this non-uniqueness via mode-filters. 

The higher weight in the initial conditions in Theorem
\ref{thstabil} vs.\,\ref{thmix1} is due to the fact that 
in Theorem \ref{thstabil} we want to extract higher order asymptotics, 
i.e., faster decay. 
The asymptotic phase-profiles in \reff{n31},\reff{n31b} and 
\reff{n31c} only depend on $k$ via $\al$ from \reff{disp-rem} 
(and on $\beta$ for \reff{n31c}). In particular, they are 
independent of the phase-speed $c_p$ and therefore are formulated 
in $x$ and $t$. 
%
% {\red Suggestion: We believe that the condition $|\phi_d|\le \eps$ 
% can be relaxed by requiring that 
% the initial interface $\phi_0$ is sufficiently flat, for instance 
% $\|\phi_0(\cdot)-\phi^*(\cdot,t_1)\|_{H^3(2)}\le \eps$ for some 
% $t_1$ sufficiently large. This will be considered elsewhere...} 
% 
\eex }
\end{Remark} 

\begin{Remark}\label{w-nr-m-rem2}{\rm Formally, we may as well 
describe the diffusive mixing of wave numbers in case $\omega''=0$, see Remark \ref{w-nr-m-rem}. However, 
then the rigorous separation of the (then unbounded) 
phase, see \S\ref{psep-sec}, becomes more difficult. 
Therefore, we will not consider this case here.
\eex}
\end{Remark}

%%%%%%%%%%%%%%%%%%%%%%%%%%%%%%%%%%%%%%%%%%%%%%%%%%%%%%%%%%%%%%%%%%%%%%%%%%%%

\subsection{The idea}\label{for-sec}
% Motivated by the asymptotic behavior of solutions to the 
% Burgers equation \reff{be}, i.e., $\pa_T q=\al\pa_X^2 q+\beta\pa_X(q^2)$, 
% which in suitable coordinates is the amplitude 
% equation for perturbations of spectrally stable wave-trains 
% $u_0(\th; k_0)$, 

The translation invariance of \reff{rd1} and the fact that by assumption 
we have periodic wave trains $u_0(\th;k)$ for wave numbers $k$ in 
a whole interval $(k_l,k_r)$ suggest to consider initial 
conditions for \reff{rd1} of the form 
\begin{gather}
|\uic(x)-u_0(k_{\pm}x+\phi_{\pm};k_{\pm})|\ra 0\quad\text{as}\quad 
x\ra\pm\infty.
\label{icform}
\end{gather}
The behavior of the corresponding solutions can be 
discussed formally if we assume that the initial phase shift $\phi_+-\phi_-$ or the initial 
wave number shift $q_+-q_-$ happens on a long spatial scale. 
We make the ansatz 
\begin{equation}\label{be-ans}
u(x,t) = u_0(k_0x-\om_0t+\Phi(X,T);k_0+\delta\pa_X\Phi(X,T))
\end{equation}
where $0<\delta\ll 1$ is a small perturbation parameter
that determines the length scale over which the wave number is
modulated by the function $\pa_X\Phi$, and where $X$ and $T$ 
are long spatial and temporal scales. Plugging \reff{be-ans} 
into \reff{rd1} and comparing equal powers in $\del$ it turns out 
(see \cite[\S4.3]{DSSS05}) that 
\begin{equation}\label{be-coord}
(X,T) = \left(\delta(x-c_\mathrm{g}t),\delta^2t\right),
\quad\text{where}\quad 
c_\mathrm{g} = \omnl^\prime(k), 
\end{equation}
are the correct spatial and temporal scales, 
and that  $q(X,T):=\pa_X\Phi(X,T)$ should satisfy  the 
Burgers equation 
\begin{equation}\label{be}
\pa_T q= \al \pa_{X}^2 q+\beta\pa_X (q^2), \quad 
\al=-\frac{1}{2} \lambda_1^{\prime\prime}(0), \ 
\beta=- \frac{1}{2} \omnl^{\prime\prime}(k), 
\end{equation}
while the phase $\Phi(X,T)$ itself satisfies the integrated Burgers
equation
\begin{equation}\label{be-phase}
\pa_T \Phi= \al \pa_{X}^2 \Phi+\beta (\pa_X\Phi)^2. 
\end{equation}
Note again that in case $\beta=0$ (Theorem \ref{thmix1} (i)) 
\reff{be} resp.~\reff{be-phase} are the linear diffusion equation.

Two questions arise: (a) What do we (formally) learn from \reff{be} resp.\,\reff{be-phase}? (b) In what sense, i.e.~in what spaces and over what timescales, do solutions of \reff{be} via \reff{be-ans} approximate solutions of \reff{rd1}, and can we give rigorous proofs for that?

To answer (a) we briefly review some well-known results about dynamics and stability in the Burgers equation in the following section. With this in mind we turn to (b). 
%In \cite{DSSS05}, (b) has been answered in the sense that \reff{be} gives an accurate description of the dynamics of \reff{rd1} on long but finite intervals. Here we consider (b) in a somewhat different direction: in a nutshell, we have a self-similar decay of the wave number perturbation $q$ in \reff{be}. Therefore we may expect some diffusive mixing of the phases in \reff{rd1}, and this is what we make precise and prove rigorously in the present paper.
One way to translate the formal analysis into rigorous results is to give 
estimates for the difference between the formal approximation
\[
U_\mathrm{approx}(x,t) =
u_0(\th+\Phi(X,T;\delta);
k+\delta\pa_X\Phi(X,T;\delta))
\]
and a true solution $u(x,t)$ of \reff{rd1}, on sufficiently long time scales. In \cite{DSSS05} this has been achieved for a variety of cases using a separation of the critical mode (the phase mode) from the exponentially damped remaining modes by Bloch wave analysis. Here, for special initial data we obtain the diffusive stability results and the mixing results from Theorems \ref{thstabil} and 
\ref{thmix1}. The proofs heavily rely on the coordinates from \cite{DSSS05}, 
which are introduced in \S\ref{psep-sec}.

%%%%%%%%%%%%%%%%%%%%%%%%%%%%%%%%%%%%%%%%%%%%%%%%%%%%%%%%%%%%%%%%%%%%%%%%%%%%

\subsection{Dynamics in the perturbed Burgers equation}\label{buri-sec}
In the spectrally stable case the amplitude equation for 
long wave modulations of the local wave number is given by the Burgers 
equation. For given classes of initial conditions 
the behavior of solutions of the Burgers equation is well understood, 
and, moreover, this behavior is stable under perturbations of the 
Burgers equation. 
Thus, before proving our results for \reff{rd1} with initial data 
\reff{icform} we briefly review some well-known results 
about the (perturbed) Burgers equation, cf.\,\cite{bkl94, msu01} and 
\cite[Sec.~3]{ue04nsb}.  This also motivates the ideas and methods of the 
proof. To keep track of $\al$ and $\beta$ 
we do not rescale \reff{be} to the standard form 
$\pa_\tau q=\pa_\xi^2 q+\pa_\xi(q^2)$. 

The Burgers equation \reff{be} has Galilean invariance: if $q$ solves 
\reff{be}, then $v=q+c$ solves 
$\pa_T v=\al\pa_X^2 v-2c\beta\pa_X v+\beta\pa_X(v)^2$ 
which can be transformed back to 
\reff{be} via $X\mapsto X+2c\beta T$. Thus, concerning 
the stability of constant solutions of \reff{be} we can 
restrict to $q\equiv 0$. 

\medskip
We add a higher order perturbation in the form of a total derivative
to \reff{be} and for notational convenience 
we take initial conditions at time $T=1$. 
Thus we consider 
\begin{equation}\label{be2}
\pa_T q = \al\pa_{X}^2 q + \beta\pa_X (q^2)
+\ga \pa_X h(q,\pa_Xq),\quad q|_{T=1}=q_0, 
\end{equation}
where for simplicity $h(a,b)=a^{d_1}b^{d_2}$ is a monomial. 
For $\ga=0$ we again have the Burgers equation. The 
perturbation is assumed to be %quasi-linear and 
of higher order. To make this 
precise we define the degree 
\begin{gather}
d_h=d_1+2d_2-3,\quad\text{
and assume that}\quad d_2\leq 1\quad\text{and}\quad d_h\geq 0. 
\label{deg-ass}
\end{gather}
%This is explained in Remark \ref{in-rem} below. 

The mean $\int_\R q(X,T)\,\rmd X$ is conserved also by the 
perturbed Burgers equation (\ref{be2}). 
Diffusive stability of $q=0$ in \reff{be2} is based on the fact that 
solutions to the linear diffusion equation in Fourier space concentrate at 
wave number $\kap=0$. Roughly speaking, 
for initial data in $L^1(\R)$ that decay like $|X|^{-n}$, 
the solutions of $\pa_T q=\al\pa_X^2 q$ fulfill 
\begin{gather}
q(X,T) = \sum^{n-1}_{j=0} T^{-(j+1)/2} {\hat q}^{(j)}_0(0) H_j 
(X/{\sqrt{T}}) + O(T^{-n/2})\quad\text{for}\quad T \to \infty,
\label{heata}
\end{gather}
where $H_j$ is a multiple of the (scaled) $j^{{\rm th}}$ Hermite function 
$H(x)=(-1)^j\pa_x^j\exp(-x^2/(4\al))$. Thus, if 
$\qhat_0(0)=\frac 1{2\pi}\int_\R q_0(X)\,\rmd X\neq 0$ then 
$\|q(\cdot,T)\|_{L^\infty} \leq C T^{-1/2} \|q_0\|_{L^1}$,  
while for  $\qhat(0)=0$ we have 
$\|q(\cdot,T)\|_{L^\infty} \leq C T^{-1} \|q_0\|_{L^1}$. 
In the second case it turns out 
that solutions to the nonlinear equation (\ref{be2}) with zero mean
have the same asymptotics as solutions to the linearization with
zero mean. Thus, both nonlinear terms
$\beta\pa_X(q^2)$ and $\ga\pa_Xh(q,\pa_Xq)$ are called 
asymptotically irrelevant. 
 
For  $\qhat(0)\neq 0$ only $\ga\pa_Xh(q,\pa_Xq)$ 
is irrelevant, and there is a nonlinear correction to the dynamics 
for \reff{be2} compared to \reff{heata}. To derive this we use 
the Cole--Hopf transformation
\[
Q(X,T) = \exp\left(\frac \beta \al \int_{-\infty}^{\sqrt{\al}X} 
q(Y,T)\,\rmd Y\right),\qquad 
q(X,T) = \frac{\sqrt{a}}\beta\frac{\pa_Y Q(Y,T)}{Q(Y,T)},\quad 
Y=X/\sqrt{\al}, 
\]
which transforms (\ref{be2}) with $\ga=0$ into the linear heat equation
$\pa_T Q = \pa_{X}^2 Q$, $Q|_{T=1}=Q_0$, 
with  $\lim_{X\to-\infty} Q_0(X) = 1$ and 
$\lim_{X\to\infty} Q_0(X) = 1+z>0$, i.e.\,
$$
\ln(1+z)=
\frac \beta \al \int_{-\infty}^{\infty} q(Y,1)\,\rmd Y.
$$
Since $\lim_{T\to\infty} Q\left(\sqrt{T}X,T\right)
= 1 + z \, \erf(X) + \rmO(1/\sqrt{T})$ 
we find that the solution $q$ to the Burgers equation (\ref{be2}) 
with $\ga=0$ satisfies 
\begin{gather}
\lim_{T\to\infty} \sqrt{T} q\left(\sqrt{T}X,T\right) =
\frac{\al}{\beta}
\frac{\rmd}{\rmd X} \ln(1+z \, \erf(X/\sqrt{\al})) =: f_z^*(X), 
\label{bua1}
\end{gather} 
with rate $\rmO(1/\sqrt{T})$. Therefore, if $\beta\neq 0$, then the
renormalized solutions converge toward a non-Gaussian limit $f_z^*(X)$. 
Again, the same behavior can be shown for \reff{be2} with 
$\ga\neq 0$. We summarize these results as follows: 

\begin{Proposition}\label{prop1}
For each $b\in(0,1/2)$, there exist $C_1,C_2,T_0>0$ such that  for 
solutions $q$ of the perturbed Burgers equation \reff{be2} the following holds. 

i) Assume that $\|q_0\|_{H^2(3)} \leq C_1$ and 
$\int_{-\infty}^\infty q_0(X)\,\rmd X = 0$. 
Then there exists a $q_{{\rm lim}}\in\R$ such that
\begin{gather}
\left\|Tq\left(\sqrt{T}X,T\right)
- q_{{\rm lim}} X \rme^{-X^2/4\al} \right\|_{H^2(3)}
\leq C_2T^{-\frac{1}{2}+b} .
\label{best1}
\end{gather}
Thus, $\|q(X,T)\|_{L^1} \leq C_2 T^{-1/2+b}$ and 
$\|q(X,T)\|_{L^\infty} \leq C_2 T^{-1+b}$. 

ii) Assume that  $A=\int_{-\infty}^\infty q_0(X)\,\rmd X \neq 0$, 
$\beta=0$, and $\|q_0\|_{H^2(2)} \leq C_1$. 
Then 
\begin{gather}
\left\|T^{1/2}q\left(\sqrt{T}X,T\right)
- \frac{A}{\sqrt{4\pi\al}}\er^{-X^2/(4\al)} \right\|_{H^2(2)}
\leq C_2T^{-\frac{1}{2}+b},  
\label{best2}
\end{gather} 
and consequently $\|q(X,T)\|_{L^\infty} \leq C_2 T^{-1/2}$. 

iii)  Assume  that $A=\int_{-\infty}^\infty q_0(X)\,\rmd X \neq 0$, 
$\beta\neq 0$, and $\|q_0\|_{H^2(2)} \leq C_1$. Then 
\begin{gather}
\left\| T^{1/2} q\left(T^{1/2}X,T\right)-f_z^*(X) \right\|_{H^2(2)}
\leq C_2T^{-1/2+b}, \text{ where } 
f_z^*(X)=\frac{\sqrt{\al}}{\beta\sqrt{4\pi}}
\frac{z\er^{-X^2/\al}}{1+z \, \erf(X/\sqrt{\al})}, 
\label{best3}
\end{gather} 
and $\ln(1+z)=\frac \beta \al \int_{-\infty}^{\infty} q_0(Y)\,\rmd Y$. 
In particular, again $\|q(X,T)\|_{L^\infty} \leq C_2T^{-1/2}$. 
\end{Proposition}

\begin{Remark}\label{smphi-rem}{\rm 
a) By translation invariance of \reff{be2}, 
we can replace $q_0$ in Proposition \ref{prop1} by $q_0(\cdot-X_0)$ for 
some $X_0\in\R$ and obtain the corresponding results for $q(X-X_0,T)$; 
w.l.o.g.~we set $X_0=0$. 
%By $T_1>T_0$ in (ii) and (iii) we assume that the initial interface between 
%$\Phi_-$ and $\Phi_+$ is sufficiently flat. 
% For instance, for $\beta=0$ the idea is that 
% for $v(X,T)=q(X,T)-Aq^*(X,T_1{+}T{-1})$ with $\int v(X,T)=0$ 
% we obtain the inhomogeneous equation 
% \begin{align}
% \pa_T v(X,T)=&\al\pa_X^2 v(X,T)+I(T,X)\notag\\
% &+\ga\pa_X\bigl[h(v(X,T)+Aq^*(X,T_1{+}T{-}1)){-}h(Aq^*(X,T_1{+}T{-}1))\bigr]
% \label{vieq}
% \end{align}
% where $I(T,X)=\ga\pa_X h(Aq^*(X,T_1+T-1))=\CO(T^{-(1+d_1/2+d_2)})$. 
% Hence, for $v$ in \reff{vieq}  
% we obtain asymptotics like for $q$ in \reff{best1}. 
%However, this requires 
%a rather detailed discussion of the higher order terms in 
%\reff{be2}, and thus also of the higher order terms in the 
%derivation of \reff{be} from 
%\reff{rd1}. Therefore we refrain from this analysis. 

b) The higher weight for $q_0$ in Proposition \ref{prop1}(i) 
compared to (ii),(iii) is due to 
the fact that we want to isolate higher order asymptotics 
(with faster decay). For this we need 
$\CF(q_0)\in H^3(2)\hookrightarrow C^2$. 

c) The profiles 
in (ii), (iii) are explicitly given in terms of $A$  
due to the conservation of $\int q\dd x$, i.e., since 
the right-hand side of \reff{be2} is a total derivative. On the other hand, 
the constant $q_{{\rm lim}}$ in (i) 
in general depends on $q_0$ in a complicated way. 

d) The local phase $\Phi$, which is related to the wave number $q$
by $q=\pa_X\Phi$, satisfies the 
(perturbed) integrated Burgers equation
\begin{equation}\label{be-phase2}
\pa_T \Phi = \al\pa_{X}^2 \Phi + \beta(\pa_X \Phi)^2
+\ga h(\pa_X\Phi,\pa_X^2\Phi), \quad \Phi(X,1)=\Phi_0(X). 
\end{equation}
For \reff{be-phase2} there exist $C_1,C_2>0$ such that we have the 
following asymptotics. 

i) If $\|\Phi_0\|_{H^3(3)}\le C_1$ then there exists an $\phi_{{\rm lim}}
=-2\al q_{{\rm lim}}\in\R$ such that
\begin{gather}\label{phik0}
\left\|T^{1/2} \Phi \left(\sqrt{T}X,T\right)
- \phi_{{\rm lim}} \rme^{-X^2/4\al} \right\|_{H^3(3)}
\leq C_2T^{-{1}/{2}+b}. 
\end{gather}
Thus, the renormalized phase converges toward a Gaussian. 

ii) If $\beta=0$ and $\Phi_0(X)\ra\Phi_\pm$ as $X\ra\pm\infty$ with 
$|\Phi_+-\Phi_-|\le C_1$ and $\|\Phi_0'\|_{H^2(2)}\leq C_1$, then 
$$
\|\Phi(\sqrt{T}X,T) - \Phi^*(X)\|_{H^3(2)}
\leq C_2T^{-1/2+b}$$  
where 
$\Phi^*(X)=\Phi_-+(\Phi_+-\Phi_-)\erf(X/\sqrt{\al})$. 

iii) If $\beta\neq 0$ and 
and $\Phi_0(X)\ra\Phi_\pm$ as $X\ra\pm\infty$ with 
$|\Phi_+-\Phi_-|\le C_1$ and $\|\Phi_0'\|_{H^2(2)}\leq C_1$, 
then 
$$\|\Phi(\sqrt{T}X,T)-\Phi^*_z(X)\|_{H^3(2)}\leq C_2T^{-1/2+b}$$ 
where $\Phi^*_z(X)=\Phi_-+\frac\al\beta\ln(1+z\erf(X/\sqrt{\al}))$, 
$\ln(1+z)=\frac\beta\al(\Phi_+-\Phi_-)$. 
\eex}
\end{Remark}

\begin{Remark}\label{in-rem}{\rm 
We briefly want to explain the reason for \reff{deg-ass} and the idea 
of (discrete) 
renormalization.  
If $\int q_0(X)\dd X\neq 0$, then, for $L>1$ chosen sufficiently large, we let  
\begin{gather}
q_n(\xi,\tau)=L^n q(L^n\xi,L^{2n}\tau). 
\label{rge1}
\end{gather}
Then $q_n$ satisfies 
\begin{gather}
  \pa_\tau q_n=\al\pa_\xi^2q_n+\beta\pa_\xi(q_n^2)+\ga L^{-n}\pa_\xi 
h_n(q_n,\pa_\xi q_n),
\label{fnd}
\end{gather}
with
\begin{gather}\label{fnd2}
h_n(q_n,\pa_\xi q_n)=L^{-(d_1+2d_2-3)n} q_n^{d_1}
(\pa_\xi q_n)^{d_2}, 
\end{gather}
where $d_h=d_1+2d_2-3\geq 0$ due to \reff{deg-ass}. 
Next, solving $\pa_T q=\al\pa_X^2 q+\beta\pa_X(v^2)+\ga\pa_X h(q,\pa_Xq)$ 
for $T\in[1,\infty)$ is equivalent to
iterating the renormalization process
\begin{gather} 
  \text{solve \reff{fnd} for } \tau\in[L^{-2},1]\text{ with initial data }
  q_n(\xi,L^{-2})=Lq_{n-1}(L\xi,1)\in\CX,
\label{r0}
\end{gather}
where $\CX$ is a suitable Banach space. Since (formally) 
$L^{-n}\pa_\xi h_n$ in \reff{fnd} 
goes to zero, in the limit $n\ra \infty$ we recover the 
linear diffusion equation (if $\beta=0$) respectively the Burgers 
equation (if $\beta\neq 0$) for $q_n$, with the known asymptotics 
\reff{heata} respectively \reff{bua1}.  Similarly, if 
$\int q_0\dd X=0$, then we scale 
\begin{gather}
q_n(\xi,\tau)=L^{2n} q(L^n\xi,L^{2n}\tau), 
\label{rge1b}
\end{gather}
and (independent of whether $\beta$ is zero or not) end up with the linear diffusion equation  in the 
respective renormalization process. 
To make this rigorous we need a suitable Banach spaces $\CX$ 
and rigorous control of the iterative process \reff{r0}, 
and again we refer to \cite{bkl94} and \cite[sec.3]{ue04nsb} for details. 
However, two observations 
are most important. (a) In \reff{fnd2} we see that each derivative in $x$  
gives an additional $L^{-1}$ in the rescaling. 
(b) The diffusive spreading in physical space corresponds to 
concentration at $\kap=0$ in Fourier space according to 
$\CF(Lu(L\cdot))(k)=\uhat(\kap/L)$. Thus, for the linear part, 
only the parabolic shape 
of the spectrum $\lam(\kap)= -\alpha \kap^2$ of $\alpha \pa_x^2$ near  
$\kap=0$ is relevant. 
\eex}
\end{Remark}

\begin{Remark}\label{w-nr-m-rem}{\rm 
If $q_0(X)\ra q_{\pm}$ for $X\ra\pm\infty$, then 
$q(\sqrt{T}X,T)=Q_0^*(X)+\CO(T^{-1/2})$ as $T\ra\infty $
for the solutions of $q_T=\al\pa_X^2 q$, where 
$Q_0^*(X)=q_-+(q_+-q_-)\erf(X/\sqrt{\al})$.  
Thus we have diffusive mixing of the wave numbers. 
Then, for $\beta=0$ and for suitable $q_0$, we have the asymptotics
\begin{gather}
q(\sqrt{T}X,T)=Q^*(X)+\CO(T^{-1/2})
\quad\text{as}\quad T\ra\infty 
\label{dm2b}
\end{gather}
for \reff{be2}, where $|Q^*(X)-Q_0^*(X)|\leq C\er^{-X^2/4}$, i.e., 
we have essentially the same asymptotics as in the linear case, 
with a small localized 
correction, see \cite{bkl94}.  On the other hand, for $\beta\neq 0$ a 
front is created, see \cite{DSSS05}. 
However, here we do not further comment on this case since below 
we focus on diffusive mixing of phases. 
\eex}\end{Remark}

%%%%%%%%%%%%%%%%%%%%%%%%%%%%%%%%%%%%%%%%%%%%%%%%%%%%%%%%%%%%%%%%%%%%%%%%%%%%

\section{The separation of the wave numbers}\label{psep-sec}

\subsection{The ansatz}

Only 
special systems such as the cGL have an $\CS^1$-symmetry and therefore 
a natural decomposition into amplitude and phase. 
Hence, 
the first step is to extract from a general reaction-diffusion system an equation for the phase, and then 
out of this for the wave number. 
We follow the 
formal derivation made in \cite{DSSS05} which uses a multi-scale
expansion which however we cannot assume a priori. 
Thus, here we proceed as follows for the reaction-diffusion system 
\reff{rd1}. As above we change coordinates via
$\th =kx-\om t$, and obtain 
\begin{equation} \label{rd23}
\pa_t u = k^2 D \pa_\th^2 u+ \om \pa_\th u + f(u).
\end{equation}
A stationary wave train $u_0(\th;k)$ of
(\ref{rd23}) with period $2\pi$ satisfies
\begin{equation} \label{stat}
k^2 D \pa_\th^2 u_0+\om \pa_\th u_0 + f(u_0) = 0.
\end{equation}
Given a smooth phase function $\phi(\vt,t)$ we seek
solutions of the form
\begin{equation}\label{bansatz}
u(\th,t) =u_0(\vt;k(1+\pa_\vt\phi(\vt,t))) + w(\vt,t),
\end{equation}
where the phase $\phi(\vt,t)$ and the coordinates $\th$ and
$\vt$ are related by
\begin{equation}\label{bzeta-xi}
\th = \vt - \phi(\vt,t). 
\end{equation}
Roughly speaking we 
require that $\pa_\vt\phi$ is small, uniformly in $\vt$, 
and that $\phi(\vt,t)$ is close to the asymptotic profile we want to extract. 
%In detail, 
% \begin{gather}
% \phi(\vt,t)=\phi^*(\vt/\sqrt{t})+\vphi(\vt,t)\quad 
% \text{with}\quad\|\vphi(\cdot,t)\|_{H^2(2)}\leq C
% \end{gather}
% for sufficiently small $C$. 
Still, \reff{bansatz} adds an additional 
degree of freedom by introducing
$\phi$; we later add additional conditions on $\phi$ and $w$,
via mode filters, to remove this additional degree of freedom again.

\begin{Remark}\label{ia-rem}{\rm 
It might seem more natural to make the ansatz
\begin{equation}\label{bdesired}
u(\th,t) = u_0(\th+\phi(\th,t);
k(1+\pa_\th\phi(\th,t))) + w(\th,t)
\end{equation}
instead of (\ref{bansatz}). However, we need to be able to
relate the dynamics of $u(\th,t)$ back to properties of the 
wave-train $u_0(\th;k)$. Thus, we would need to express $u(\th,t)$ 
in terms of $\vt=\th+\phi(\th,t)$, i.e., 
\[
u_0(\th+\phi(\th,t);k(1+\pa_\th\phi(\th,t))) \longmapsto
u_0(\vt;k(1+\pa_\th\phi(\th(\vt,t),t)))
\]
which involves the inverse $\th(\vt,t)$ of the function
$\vt=\th+\phi(\th,t)$. The occurrence of this inverse
would have made the forthcoming analysis much more complicated.
\eex}\end{Remark}

\begin{Remark}\label{inter-rem}{\rm 
Suppose that we found a phase function $\phi(\vt,t)$ with small
derivative $\pa_\vt\phi(\vt,t)$ so that
(\ref{bansatz}) satisfies (\ref{rd23}). Using the implicit
function theorem, we can then, a posteriori, solve
(\ref{bzeta-xi}) for $\vt$ as a function of $\th$ which is
of the form $\vt=\th+\check{\phi}(\th,t)$, where
\[
\check{\phi}(\th,t)= \phi(\vt,t)= \phi(\th+\check{\phi}(\th,t),t).
\]
In particular, we see that
\begin{equation}\label{bansatz:back}
u_0(\vt;k(1+\pa_\th\phi(\vt,t)))=u_0(\th+\phi(\th+\check{\phi}(\th,t),t);
k(1+\pa_\th\phi(\th+\check{\phi}(\th,t),t)))
\end{equation}
and
\begin{align*}
\frac{\rmd}{\rmd\th}\phi(\th+\check{\phi}(\th,t),t)
&= (1+\pa_\th\check{\phi}(\th,t)) \pa_\th\phi(\th+\check{\phi}(\th,t),t)\\
&= \pa_\th\phi(\th+\check{\phi}(\th,t),t)
+\rmO(|\pa_\th\phi(\th+\check{\phi}(\th,t),t)|^2).
\end{align*}
Thus, to leading order, the solution (\ref{bansatz:back}) is 
of the desired form (\ref{bdesired}) with $\phi(\th,t)$
replaced by $\phi(\th+\check{\phi}(\th,t),t)$.\eex}\end{Remark}

We now substitute the ansatz (\ref{bansatz}) into
(\ref{rd23}) and derive the resulting PDE in $\vt$. We use
the notation
\begin{equation}\label{d:uophi}
u_0^\phi := u_0(\vt;k(1+\pa_\vt\phi)),\quad
\pa_j u_0^\phi
:= (\pa_j u_0)^\phi
:= (\pa_j u_0)(\vt;k(1+\pa_\vt\phi)), \quad
j=\vt,k. 
\end{equation}

Assuming that $\pa_\vt\phi$ is small, we obtain 
\begin{xalignat*}{2}
\frac{\rmd\vt}{\rmd\th} & = 
\frac{1}{1-\pa_\vt\phi}, &
\frac{\rmd\vt}{\rmd t} &= 
\frac{-\pa_t\phi}{1-\pa_\vt\phi}, \\
\frac{\rmd }{\rmd\th} & = 
\frac{1}{1-\pa_\vt\phi}\,\frac{\rmd}{\rmd\vt}, & 
\frac{\rmd^2}{\rmd\th^2} &= 
\left(\frac{1}{1-\pa_\vt\phi}\;\frac{\rmd}{\rmd\vt}\right)^2 , 
\end{xalignat*}
and therefore
\begin{eqnarray*}
\frac{\rmd u}{\rmd\th} & = &
\frac{1}{1-\pa_\vt\phi}\pa_\vt\,u_0^\phi
+ \frac{k\pa_\vt^2\phi}{1-\pa_\vt\phi}\, 
\pa_k u_0^\phi, \\
\frac{\rmd^2 u}{\rmd\th^2} & = &
\left(\frac{1}{1-\pa_\vt\phi}\,\frac{\rmd}{\rmd\vt}
+ \frac{k\pa_\vt^2\phi}{1-\pa_\vt\phi}\,
\frac{\rmd}{\rmd k}\right)^2 u_0^\phi, \\
\frac{\rmd u}{\rmd t} & = &
\frac{-\pa_t\phi}{1-\pa_\vt\phi}\,\pa_\vt u_0^\phi
+ k\left(
-\frac{\pa_\vt^2\phi\,\pa_t\phi}{1-\pa_\vt\phi}
+\pa_\vt\pa_t\phi\right)\pa_k u_0^\phi, 
\end{eqnarray*}
and
\begin{eqnarray*}
\frac{\rmd w}{\rmd t} & = &
\frac{\pa w}{\pa t} - \frac{\pa w}{\pa\vt}\,
\frac{\pa_t\phi}{1-\pa_\vt\phi}, \quad
\frac{\rmd w}{\rmd\th} = 
\frac{1}{1-\pa_\vt\phi}\,\frac{\pa w}{\pa\vt}, \quad
\frac{\rmd^2 w}{\rmd\th^2} = 
\left(\frac{1}{1-\pa_\vt\phi}\,
\frac{\rmd}{\rmd\vt}\right)^2 w.
\end{eqnarray*}
Thus 
\begin{eqnarray} \label{wilma}
\lefteqn{-\frac{\pa_t\phi}{1-\pa_\vt\phi}\,
\pa_\vt u_0^\phi - k\left(\pa_\vt^2\phi\,
\frac{\pa_t\phi}{1-\pa_\vt\phi}- \pa_\vt\pa_t\phi\right)\pa_k u_0^\phi
+ \pa_t w - \frac{\pa_t\phi}{1-\pa_\vt\phi}\,\pa_\vt w} 
\\ \nonumber 
& = &k^2 D \left( \left(\frac{1}{1-\pa_\vt\phi}\,\frac{\pa}{\pa\vt}
+ \frac{k\pa_\vt^2\phi}{1-\pa_\vt\phi}\,\frac{\pa}{\pa k}\right)^2 u_0^\phi
+ \left( \frac{1}{1-\pa_\vt\phi}\,\frac{\pa}{\pa\vt} \right)^2 w \right) 
\\ \nonumber & &
+\om\,\frac{1}{1-\pa_\vt\phi}\,\left(\pa_\vt u_0^\phi
+ k(\pa_\vt^2\phi)\pa_{k}u_0^\phi+ \pa_\vt w \right) \\ \nonumber & &
- (k^2 D \pa_\vt^2 u_0- \om \pa_\vt u_0 + f(u_0)) + f(u_0^\phi+w)
\end{eqnarray}
where we used (\ref{stat}) in the last equation.

Our goal is to separate the critical modes, which involve the
dynamics of $\phi$, from the damped noncritical modes using the
eigenfunctions of the linearization $\CL(k)$. 
This is done via Bloch waves which we introduce next.

%%%%%%%%%%%%%%%%%%%%%%%%%%%%%%%%%%%%%%%%%%%%%%%%%%%%%%%%%%%%%%%%%%%%%%%%%%%%

\subsection{Bloch wave analysis}\label{bw-sec}

Bloch wave transform $\mathcal{J}$ is 
a generalization of Fourier transform $\mathcal{F}$. We briefly 
review the main properties and refer to \cite{RS72,gs98,scar99,DSSS05} 
for proofs and further details. From now on we use a slightly
rescaled Fourier transform, namely 
\begin{gather}
\hat{w}(\ell)=({\cal F}w)(\ell)=\frac{1}{2\pi k}\int_{-\infty}^\infty 
\er^{-\ri \ell \vt/k}w(\vt)\dd\vt,
\quad w(\vt)=({\cal F}^{-1}\hat{w})(\vt)=
\int_{-\infty}^\infty \er^{\ri \ell\vt/k}\hat{w}(\ell)\dd\ell, 
\label{sft}
\end{gather}
and thus, denoting the classical Fourier transform (where $k=1$) 
by ${\cal F}_1$, 
\huga{
({\cal F}w)(\ell)=\frac 1 k ({\cal F}_1w)(\ell/k)\quad \text{ and }\quad 
({\cal F}\hat w)(\vt)=({\cal F}_1^{-1}\hat{w})(\vt/k).
\label{sft2}
}
Then, for  sufficiently smooth and rapidly enough decaying functions
$w$, we have
\begin{align}
(\CJ^{-1}&\wti)(\vt):=w(\vt)=
\int_{-\infty}^\infty\rme^{\rmi\ell\vt/k}\hat{w}(\ell)\,\rmd\ell
=\sum_{j\in\Z}\int_{-k/2}^{k/2}\rme^{\rmi\vt(\ell+jk)/k}
      \hat{w}(\ell+jk)\,\rmd\ell \notag \\
= & \int_{-k/2}^{k/2} \rme^{\rmi\ell\vt/k}
      \left[\sum_{j\in\Z}\rme^{\rmi j\vt}\hat{w}(\ell+jk)\right]
      \,\rmd\ell=
\int_{-k/2}^{k/2}\rme^{\rmi\ell\vt/k}
       \tilde{w}(\vt,\ell)\,\rmd\ell \label{j-def}
\end{align}
where 
\begin{gather}
\wti(\vt,\ell)=(\CJ w)(\vt,\ell)=
\sum_{j\in\Z}\rme^{\rmi j\vt}\hat{w}(\ell+jk). 
\label{jdef}
\end{gather}
Similar to the Fourier
transform, the Bloch transform can be defined for tempered
distributions. By construction, 
\begin{gather}
\wti(\vt+2\pi,\ell)=\wti(\vt,\ell)\quad\text{and}\quad
\wti(\vt,\ell+k)=\er^{\ri\vt}\wti(\vt,\ell),
\label{bper}
\end{gather}
such that we can restrict ourselves to $\ell\in[-k/2,k/2)$. 
The Bloch transform of the product of two functions $w_1$ and
$w_2$ in $\vt$-space is given by the convolution
\begin{equation}\label{*}
\mathcal{J}[w_1 \cdot w_2](\vt,\ell)
= [\tilde{w}_1 * \tilde{w}_2](\vt,\ell)
= \int_{-k/2}^{k/2} \tilde{w}_1(\vt,\ell-\tilde{\ell})
  \tilde{w}_2(\vt,\tilde{\ell})\,\rmd\tilde{\ell}
\end{equation}
of their Bloch transforms $\tilde{w}_1$ and $\tilde{w}_2$ in Bloch space, 
where \reff{bper} is used for $|\ell-\tilde{\ell}|>k/2$. 
The analytic properties of the Bloch transform are based on a
generalization of Parseval's identity
\[
\int_{-\infty}^\infty |u(\vt)|^2\,\rmd\vt =
2\pi k \int_0^{2\pi} \int_{-k/2}^{k/2} |\tilde{u}(\vt,\ell)|^2
\,\rmd\ell\,\rmd\vt.
\]
As a consequence, Bloch wave transform is an isomorphism between 
$H^{m_2}(m_1)$, 
and the space $B^{m_1}(m_2)$ of functions
$\tilde{u}(\vt,\ell)$ that are $2\pi$-periodic w.r.t.\ $\vt$,
satisfy (\ref{bper}), and whose norm
\[
\|\tilde{u}\|_{B^{m_1}(m_2)} =
\sum_{j=0}^{m_1} \sum_{i=0}^{m_2} \int_0^{2\pi} \int_{-k/2}^{k/2}
|\pa_\ell^j \pa_\vt^i \tilde{u}(\vt,\ell)|^2
\,\rmd\ell\,\rmd\vt
\]
is finite. 
We now collect a few more properties; see, e.g., \cite{gs98} or 
\cite[\S5.2]{DSSS05}  for more details and proofs. 
\begin{Remark}\label{btrem}{\rm 
a) If $w_1(\vt)$ is $2\pi$-periodic in $\vt$ and the
support of the Fourier transform $\hat{w}_2$ of a complex-valued function
$w_2(\vt)$ lies in $(-1/2,1/2)$, then we have 
\begin{equation}\label{fbconv}
\mathcal{J}[w_1 w_2](\vt,\ell) = w_1(\vt) \hat{w}_2(\ell).
\end{equation} 
Due to \reff{fbconv} Bloch transform is useful to analyze differential 
operators with
spatially-periodic coefficients, which in Bloch space become 
multiplication operators. 

b) Since we are interested in functions which do not necessarily
decay to zero at infinity, we employ a method already used 
in \cite{Schn94ZAMP} to extend
multiplication operators from the space $L^2$ of square-integrable functions to
the space $L^2_\mathrm{ul}$ of uniformly locally square-integrable functions
equipped with the norm
$\ds \|u\|_{L^2_\mathrm{ul}} =
\sup_{x\in\R} \int_{x}^{x+1} |u(y)|^2 \,\rmd y.$
We recall that
\[
H^m_\mathrm{ul} =
\{ u:\R\rightarrow\R;\;
\|u\|_{H^m_\mathrm{ul}} = \|u\|_{H^m(x,x+1)} < \infty
\mbox{ with }
\lim_{y\rightarrow0} \|u-T_y u\|_{H^m_\mathrm{ul}} \rightarrow0 \}
\]
where $[T_y u](x)=u(x+y)$. Now let $m,s\in\Z$ with $m+s\geq0$ and $m\geq0$, 
and consider a function
\[
\tilde{\mathcal{M}}: \quad
\R \longrightarrow
L(H^{m+s}_\mathrm{per}(0,2\pi),H^m_\mathrm{per}(0,2\pi)), \qquad
\ell \longmapsto \tilde{\mathcal{M}}(\ell)
\]
which is $\mathcal{C}^2$ in the Bloch wave
number $\ell$. Then $\tilde{\mathcal{M}}$ defines a bounded
operator $\mathcal{M}:H^{m+s}_\mathrm{ul}\rightarrow H^m_\mathrm{ul}$
with 
\huga{\label{opiest}
\|\mathcal{M}\|_{L(H^{m+s}_\mathrm{ul},H^m_\mathrm{ul})} \leq C(m,s)
\|\tilde{\mathcal{M}}\|_{\mathcal{C}^2_\mathrm{b}%
((-k/2,k/2),L(H^{m+s}_\mathrm{per},H^m_\mathrm{per}))}.
}
Clearly, this can  be extended to multi-linear operators.\eex}
\end{Remark}
%%%%%%%%%%%%%%%%%%%%%%%%%%%%%%%%%%%%%%%%%%%%%%%%%%%%%%%%%%%%%%%%%%%%%%%%%%%%

\subsection{Mode filters, and separation into critical and noncritical modes}
\label{s:sep}
Our goal is to separate the dynamics of the eigenmodes
$\tilde{v}_1(\vt,\ell)$ associated with the critical eigenvalues
$\lambda_1(\ell)$ from the remaining modes, which are linearly 
exponentially damped and therefore called non-critical. 
We use mode filters to obtain this splitting.

Due to Hypotheses \ref{hyp1} there exists a number $\ell_1$ with
$0<\ell_1\ll1$ so that the eigenvalue $\lambda_1(\ell)$ of
$\tilde{\CL}(\ell)$ is bounded away from the rest of the spectrum for 
$|\ell|<\ell_1$. Therefore, there exists an $\CLt(\ell)$-invariant
projection
\[
\tilde{Q}^\mathrm{c}(\ell) = \frac{1}{2\pi\rmi}
\int_\Gamma [\lambda-\tilde{\CL}(\ell)]^{-1}\,\rmd\lambda
\]
onto the space spanned by $\tilde{v}_1(\vt,\ell)$, where
$\Gamma\subset\C$ is a small circle that surrounds $\lambda_1(\ell)$
counter-clockwise in the complex plane and does not intersect the rest
of the spectrum of $\CL(\ell)$ for this fixed $\ell$. For $\ell=0$ 
we have 
$$
\tilde{Q}^\mathrm{c}(0)\vti(\cdot,0)=\langle \uad,\vti\rangle\vti_1(\cdot,0), 
$$
and similarly $\tilde{Q}^\mathrm{c}(\ell)$ can be expressed by using 
the scalar product with the adjoint $\tilde{u}_{{\rm ad}}(\cdot,\ell)$ in 
Bloch space. 

We choose a nonincreasing (for $\ell\ge 0$) 
$\mathcal{C}^\infty_0$-cutoff function
$\chi:\R\to[0,1]$ with
\begin{gather}\label{chidef}
\chi(\ell) =
\left\{ \begin{array}{lll}
1 & \mbox{for} & |\ell| \leq 1,\\
0 & \mbox{for} & |\ell| \geq 2,\\
\end{array} \right.
\end{gather}
and define
\begin{xalignat*}{2}
\tilde{P}_\mathrm{fs}^\mathrm{c}(\ell) 
&= \tilde{Q}^\mathrm{c}(\ell) \chi\left(\frac{4\ell}{\ell_1}\right), &
\tilde{P}_\mathrm{fs}^\mathrm{s}(\ell) 
&:= 1-\tilde{Q}^\mathrm{c}(\ell)\chi\left(\frac{4\ell}{\ell_1}\right),\\
\tilde{P}_\mathrm{mf}^\mathrm{c}(\ell) 
&= \tilde{Q}^\mathrm{c}(\ell) \chi\left(\frac{8\ell}{\ell_1}\right),&
\tilde{P}_\mathrm{mf}^\mathrm{s}(\ell) 
&:= 1-\tilde{Q}^\mathrm{c}(\ell)\chi\left(\frac{8\ell}{\ell_1}\right),\\
\intertext{and}
\tilde{P}^\mathrm{c}(\ell) 
&= \tilde{Q}^\mathrm{c}(\ell) \chi\left(\frac{2\ell}{\ell_1}\right),&
\tilde{P}^\mathrm{s}(\ell) 
&:= 1-\tilde{Q}^\mathrm{c}(\ell)\chi\left(\frac{16\ell}{\ell_1}\right).
\end{xalignat*}
These operators commute for each fixed $\ell$ and satisfy
\begin{gather}\label{mfprop}
\barr{l}
(1-\tilde{P}^\mathrm{c})\tilde{P}_\mathrm{fs}^\mathrm{c}
        = 0 = (1-\tilde{P}_\mathrm{fs}^\mathrm{c})
        \tilde{P}_\mathrm{mf}^\mathrm{c}, \quad
(1-\tilde{P}^\mathrm{s})\tilde{P}_\mathrm{fs}^\mathrm{s}
        = 0 = (1-\tilde{P}^\mathrm{s})\tilde{P}_\mathrm{mf}^\mathrm{s}, \\
\tilde{P}_\mathrm{fs}^\mathrm{c}+\tilde{P}_\mathrm{fs}^\mathrm{s} = 1, \quad
\tilde{P}_\mathrm{mf}^\mathrm{c}+\tilde{P}_\mathrm{mf}^\mathrm{s} = 1.
\earr
\end{gather}
We define scalar-valued operators
$\tilde{p}_\mathrm{fs}^\mathrm{c}(\ell)$ and
$\tilde{p}_\mathrm{mf}^\mathrm{c}(\ell)$ implicitly by
\begin{equation}\label{pmf}
[\tilde{p}_\mathrm{fs}^\mathrm{c}(\ell)\tilde{u}]\tilde{v}_1(\cdot,\ell)
        = \tilde{P}_\mathrm{fs}^\mathrm{c}(\ell)\tilde{u}, \qquad
[\tilde{p}_\mathrm{mf}^\mathrm{c}(\ell)\tilde{u}]\tilde{v}_1(\cdot,\ell)
        = \tilde{P}_\mathrm{mf}^\mathrm{c}(\ell)\tilde{u}. 
\end{equation}
Remark \ref{btrem}b) implies that each of the operators above extends to a
bounded operator on $H^{m+s}_\mathrm{ul}$. The resulting operators will be
denoted by the same letter but with the superscript $\tilde{~}$ being dropped.

The mode filters $p_\mathrm{mf}^\mathrm{c}$ and $P_\mathrm{mf}^\mathrm{s}$
are now used to separate the critical and
noncritical modes in (\ref{wilma}), while 
$p_\mathrm{fs}^\mathrm{c}$ and $P_\mathrm{fs}^\mathrm{s}$ 
are used to limit the Fourier support of the critical modes.
We write (\ref{wilma}) given by
\begin{eqnarray}
\lefteqn{-\frac{\pa_t\phi}{1-\pa_\vt\phi}\,
\pa_\vt u_0^\phi
- k\left(\pa_\vt^2\phi\,
\frac{\pa_t\phi}{1-\pa_\vt\phi}
- \pa_\vt\pa_t\phi\right)\pa_k u_0^\phi
+ \pa_t w
- \frac{\pa_t\phi}{1-\pa_\vt\phi}\,\pa_\vt w}
\nonumber \\ \label{rd-pw} & = &
k^2 D \left( \left(
\frac{1}{1-\pa_\vt\phi}\,\frac{\pa}{\pa\vt}
+ \frac{k\pa_\vt^2\phi}{1-\pa_\vt\phi}\,
\frac{\pa}{\pa k}\right)^2 u_0^\phi
+ \left( \frac{1}{1-\pa_\vt\phi}\,
\frac{\pa}{\pa\vt} \right)^2 w \right) \\ \nonumber & &
- \om\,\frac{1}{1-\pa_\vt\phi}\,
\left(\pa_\vt u_0^\phi
+ k(\pa_\vt^2\phi)\pa_{k}u_0^\phi
+ \pa_\vt w \right) \\ \nonumber & &
- (k^2 D \pa_\vt^2 u_0
- \om \pa_\vt u_0 + f(u_0)) + f(u_0^\phi+w),
\end{eqnarray}
as
\begin{equation}\label{rd-befmf}
[-B_0 + B_1(\pa_\vt\phi,w)] \pa_t\phi + \pa_t w
= -\CL_i \pa_\vt\phi + \CL w
+ G(\pa_\vt\phi,w),
\end{equation}
where
\begin{eqnarray}
B_0 \pa_t\phi & = &
  (\pa_\vt u_0 - k \pa_k u_0\pa_\vt)\pa_t\phi, 
\nonumber \\ \nonumber
\CL_i \pa_\vt\phi & = &
  - \CL(k\pa_\vt\phi\,\pa_{k}u_0)
  + k^2 D (2\pa_\vt\phi\,\pa_\vt^2 u_0
  + \pa_\vt^2\phi\,\pa_\vt u_0)
  - \om\pa_\vt\phi\,\pa_\vt u_0
\\ \label{l0tilde} & = &
  k\left[ \CL(\ri\pa_\vt\phi\,\pa_\ell v_1)
  + k D (2\pa_\vt\phi\,\pa_\vt^2 u_0
  + \pa_\vt^2\phi\,\pa_\vt u_0)
  - c_\mathrm{p}\pa_\vt\phi\,\pa_\vt u_0 \right], 
\\ \nonumber
B_1 (\pa_\vt\phi,w) \pa_t\phi & = &
  \left(\pa_\vt u_0
  -\frac{\pa_\vt u_0^\phi}{1-\pa_\vt\phi}\right)
    \pa_t\phi
  - k\left(\frac{\pa_\vt^2\phi\,\pa_k u_0^\phi}%
    {1-\pa_\vt\phi}\pa_t\phi
  + (\pa_k u_0-\pa_k u_0^\phi)\pa_\vt\pa_t\phi\right)%\notag\\
 - \frac{\pa_\vt w}{1-\pa_\vt\phi}\,\pa_t\phi,\notag
\end{eqnarray}
and where $G$ is contains the remaining terms. 
In the calculation above,
we used that $\pa_\ell v_1=\ri\pa_k u_0$, see \reff{dlv}. 
The symbol 
$\CL_i$ is used since in the critical modes $\pa_\vt\CL_i$ corresponds 
to $\CL$, see \reff{dcheck} below, i.e., $\CL_i$ resembles 
an integration of $\CL$. Clearly, 
$$
B_1(\pa_\vt\phi,w)  = 
   \rmO(|\pa_\vt\phi|+|w|), \quad 
G(\pa_\vt\phi,w)  = \rmO(|\pa_\vt\phi|^2+|w|^2).
$$

Our goal is to replace (\ref{rd-befmf}) with the system
\begin{eqnarray}\label{rd-mf}
\pa_t P_\mathrm{fs}^\mathrm{c} B_0 \phi & = &
  P_\mathrm{fs}^\mathrm{c} \CL_i \pa_\vt\phi
+ P_\mathrm{mf}^\mathrm{c} B_1(\pa_\vt\phi,w) \pa_t\phi
- P_\mathrm{mf}^\mathrm{c} G(\pa_\vt\phi,w), 
\\ \nonumber
\pa_t w & = & \CL w
+ P_\mathrm{fs}^\mathrm{s} B_0 \pa_t\phi
- P_\mathrm{fs}^\mathrm{s} \CL_i \pa_\vt\phi
- P_\mathrm{mf}^\mathrm{s} B_1(\pa_\vt\phi,w) \pa_t\phi
+ P_\mathrm{mf}^\mathrm{s} G(\pa_\vt\phi,w)
\end{eqnarray}
for $(\phi,w)$. Subtracting the first from the second equation and using
(\ref{mfprop}), we see that solutions of (\ref{rd-mf}) give solutions of
(\ref{rd-befmf}). Alternatively, we may consider the system
\begin{eqnarray}\label{rd-mfp}
\pa_t p_\mathrm{fs}^\mathrm{c} B_0 \phi & = &
  p_\mathrm{fs}^\mathrm{c} \CL_i \pa_\vt\phi
+ p_\mathrm{mf}^\mathrm{c} B_1(\pa_\vt\phi,w) \pa_t\phi
- p_\mathrm{mf}^\mathrm{c} G(\pa_\vt\phi,w), 
\\ \nonumber
\pa_t w & = & \CL w
+ P_\mathrm{fs}^\mathrm{s} B_0 \pa_t\phi
- P_\mathrm{fs}^\mathrm{s} \CL_i \pa_\vt\phi
- P_\mathrm{mf}^\mathrm{s} B_1(\pa_\vt\phi,w) \pa_t\phi
+ P_\mathrm{mf}^\mathrm{s} G(\pa_\vt\phi,w)
\end{eqnarray}
for $(\phi,w)$, where the first equation is now scalar-valued. Inspecting
(\ref{pmf}) 
we see that
(\ref{rd-mf}) and (\ref{rd-mfp}) are equivalent. 
We shall require that $(\phi,w)$ satisfy
\begin{equation}\label{rd-cond-phi}
\mathrm{supp}\,\mathcal{F}[\phi] \subset \mathcal{I} :=
\left\{\ell;\; \chi\left(4\ell/\ell_1\right)=1 \right\}
\end{equation}
and
\begin{equation}\label{rd-cond-w}
(1-P^\mathrm{s}) w = 0
\end{equation}
for all $t\geq 1$. Since $P^\mathrm{s}$ commutes with $\CL$, it follows
from (\ref{mfprop}) and (\ref{rd-mfp}) that (\ref{rd-cond-w}) holds
for all $t>1$ if it is true for $t=1$.

It remains to check whether (\ref{rd-cond-phi}) is respected by
(\ref{rd-mfp}) and to calculate the operator $p_\mathrm{fs}^\mathrm{c}B_0$
to see whether (\ref{rd-mfp}) is a proper evolution equation. 
Due to the properties of the
multiplier $p_\mathrm{mf}^\mathrm{c}$, we know that
\[
\mathrm{supp}\,\mathcal{F}\left[p_\mathrm{mf}^\mathrm{c}
(B_1(\pa_\vt\phi,w)\pa_t\phi-G(\pa_\vt\phi,w))
\right] \Subset \mathcal{I}
\]
for any sufficiently smooth function $\phi$. 
>From (\ref{l0tilde}) we find that the operators $B_0$ and 
$\CT_i$ have $2\pi$-periodic coefficients 
in $\vt$ and are multipliers in Bloch space which allows us
to use Remark \ref{btrem}. For any function $\phi$ that satisfies
(\ref{rd-cond-phi}), we then obtain
\begin{align*}
\tilde{P}_\mathrm{fs}^\mathrm{c} \mathcal{J}\left[B_0\phi\right] &= 
        \tilde{P}_\mathrm{fs}^\mathrm{c}(\ell)\mathcal{J}[B_0\phi](\vt,\ell)
 \stackrel{(\ref{fbconv})}{=} 
        \hat{\phi}(\ell) \chi\left({4\ell}/{\ell_1}\right)
        \tilde{Q}^\mathrm{c}(\ell)
        \left(\pa_\vt u_0(\vt) + \rmO(\ell)\right)\\
&=   \hat{\phi}(\ell) \chi\left({4\ell}/{\ell_1}\right)
        \left(1+\rmO(\ell)\right) \tilde{v}_1(\vt,\ell)
\stackrel{(\ref{rd-cond-phi})}{=} 
        \left[ \left(1+\rmO(\ell_1)\right) \hat{\phi}(\ell) \right]
        \tilde{v}_1(\vt,\ell),
\end{align*}
where the $\rmO(\ell_1)$-term is a multiplier and 
$\left[\left(1+\rmO(\ell_1)\right) \hat{\phi}\right]$ has support 
in $\mathcal{I}$. Therefore, using the definition
(\ref{pmf}) of $p_\mathrm{fs}^\mathrm{c}$ and denoting the operator
associated with the $\rmO(\ell_1)$-term by $B_2$, we get
\begin{equation}\label{b0exp}
p_\mathrm{fs}^\mathrm{c} B_0 \phi = (1+B_2) \phi
\end{equation}
for all $\phi$ that satisfy (\ref{rd-cond-phi}), where $B_2$ has norm 
$\|B_2\|=\rmO(\ell_1)$ and respects (\ref{rd-cond-phi}), i.e.\
$\mathrm{supp}\,\mathcal{F}[B_2\phi]\subset\mathcal{I}$. Since similar
arguments apply to the multiplier $\CL_i$, 
(\ref{rd-cond-phi}) is indeed preserved by (\ref{rd-mfp}). 

For all $(\phi,w)$ for which $(\pa_\vt\phi,w)$ is small
and $\phi$ satisfies (\ref{rd-cond-phi}), the first equation of
(\ref{rd-mfp}) can be written as
\[
\pa_t \phi =
\left[1+B_2+p_\mathrm{mf}^\mathrm{c}B_1(\pa_\vt\phi,w)\right]^{-1}
\left[ p_\mathrm{fs}^\mathrm{c} \CL_i \pa_\vt\phi
- p_\mathrm{mf}^\mathrm{c} G(\pa_\vt\phi,w) \right].
\]
Substituting this expression for $\pa_t\phi$ into the second equation
of (\ref{rd-mfp}) for $w$, we arrive at the system
\begin{align}
\pa_t \phi = &
\left[1+B_2+p_\mathrm{mf}^\mathrm{c}B_1(\pa_\vt\phi,w)\right]^{-1}
\left[ p_\mathrm{fs}^\mathrm{c} \CL_i \pa_\vt\phi
- p_\mathrm{mf}^\mathrm{c} G(\pa_\vt\phi,w) \right], \label{mfphi}\\
\pa_t w = & \CL w
- P_\mathrm{fs}^\mathrm{s} \CL_i \pa_\vt\phi
+ P_\mathrm{mf}^\mathrm{s} G(\pa_\vt\phi,w)\notag\\ 
&+ \left[ P_\mathrm{fs}^\mathrm{s} B_0
        {-}P_\mathrm{mf}^\mathrm{s} B_1(\pa_\vt\phi,w) \right]
\left[1{+}B_2{+}p_\mathrm{mf}^\mathrm{c}B_1(\pa_\vt\phi,w)\right]^{-1}
\left[ p_\mathrm{fs}^\mathrm{c} \CL_i \pa_\vt\phi
{-}p_\mathrm{mf}^\mathrm{c} G(\pa_\vt\phi,w) \right].\label{we1}
\end{align}
Thus we have a splitting of the critical modes 
$\phi$ and the noncritical modes $w$.

%%%%%%%%%%%%%%%%%%%%%%%%%%%%%%%%%%%%%%%%%%%%%%%%%%%%%%%%%%%%%%%%%%%%%%%%%%%%

\subsection{The system for wave numbers and damped modes}
\label{wnr-sec}
We now replace $\phi$ by $\psi=\pa_\vt\phi$ and obtain
\begin{align}
\pa_t \psi = & \pa_\vt
\left[1+B_2+p_\mathrm{mf}^\mathrm{c}B_1(\psi,w)\right]^{-1}
\left[ p_\mathrm{fs}^\mathrm{c} \CL_i \psi
- p_\mathrm{mf}^\mathrm{c} G(\psi,w) \right], \label{skew1}\\
\pa_t w = & \CL w- P_\mathrm{fs}^\mathrm{s} \CL_i \psi
+ P_\mathrm{mf}^\mathrm{s} G(\psi,w)\notag\\ 
&+ \left[ P_\mathrm{fs}^\mathrm{s} B_0
        {-}P_\mathrm{mf}^\mathrm{s} B_1(\psi,w) \right]
\left[1{+}B_2{+}p_\mathrm{mf}^\mathrm{c}B_1(\psi,w)\right]^{-1}
\left[ p_\mathrm{fs}^\mathrm{c} \CL_i \psi
{-}p_\mathrm{mf}^\mathrm{c} G(\psi,w) \right], \label{skew2}
\end{align}
which we also write in short as
\begin{equation} \label{ws1}
\pa_t \mathcal{V} = \Lambda \mathcal{V} + F(\mathcal{V}),
\end{equation}
where $\mathcal{V}=(\psi,w)$, $\Lambda$ is a linear operator, and
$F(\mathcal{V})=\rmO(|\mathcal{V}|^2)$. 
We now prove that the spectrum of the operator
\[
\pa_\vt (1+B_2)^{-1} p^\mathrm{c}_\mathrm{mf} \CL_i
\]
near $\lambda=0$ is approximately given by the linear dispersion curve
$\lam_1(\ell)$ with the associated eigenmodes given
approximately by the Fourier modes $\exp(-\rmi\ell\vt/k)$. This follows 
from 
\begin{align*} 
\CL_i(\er^{\ri\ell\vt/k})&= k\bigl[\CL(\er^{\ri \ell\vt/k}\ri\pa_\ell\vti_1)
+(kD(2\pa_\vt^2u_0+\ri(\ell/k)\pa_\vt u_0)
-\ri c_p\pa_\vt u_0)\er^{\ri\ell\vt/k}\bigr], 
\end{align*}
and therefore 
$\ds (\tilde{p}_\mathrm{mf}^\mathrm{c}\CJ\CL_i(\er^{\ri\ell\vt/k}))(\ell)
=\chi(\frac{8\ell}\ell_1)\ri k
\bigl[\lam_1'(0)-\ri\ell(2kD\pa_\vt\pa_k\vti_1+\pa_\vt u_0)+\CO(\ell^2)\bigr]
\er^{\ri\ell\vt/k}.$
Since $1+B_2(\ell)=1+\CO(\ell)$ as a multiplier and 
$\pa_\vt\er^{\ri\ell\vt/k}=\ri(\ell/k)\er^{\ri\ell\vt/k}$ we find 
\begin{gather}
\CJ(\pa_\vt(1+B_2)^{-1}p_\mathrm{mf}^\mathrm{c}\CL_i\er^{\ri\ell\vt/k})(\ell)
=\chi(\frac{8\ell}{\ell_1})(\lam_1(\ell)+\CO(\ell^3))\er^{\ri\ell\vt/k}.
\label{dcheck}
\end{gather}

% The diffusive variable $\psi$ appears linearly in \reff{skew2} but 
% only with support of $\psiti$ in $\{|\ell|>\ell_1/8\}$, hence 
% it is in fact linearly exponentially damped in \reff{skew2}. Nevertheless, 
For notational convenience we diagonalize the linear 
part of \reff{skew1},\reff{skew2} by setting 
\begin{gather}
\bpm v^c\\ v^s\epm=\CS^{-1}\bpm \psi\\ w\epm =\bpm 1&0\\ -\CS_1 &1\epm 
\bpm \psi\\ w\epm, \label{sdef}
\end{gather}
where $\tilde{\CS}_1\in C^\infty([-k/2,k/2),L(\C,H^{m}(\CT_{2\pi})))$ 
is a multiplier with 
$\supp\tilde{\CS}_1\subset\{\ell_1/8<|\ell|<\ell_1/4\}$. 
Thus, 
$v^\mathrm{c}=\psi$ and $P^\mathrm{s}v^\mathrm{s}=v^\mathrm{s}$, 
and, by definition, 
\begin{equation}\label{s}
\CS^{-1} \Lambda \CS =
\mathrm{diag}(\lambda^\mathrm{c},\Lambda^\mathrm{s}), 
\end{equation}
with $\lam^c(\ell)=\chi(\frac{8\ell}{\ell_1})(\lam_1(\ell)+\CO(\ell^3))$, 
cf.~\reff{dcheck}. 
In these coordinates, (\ref{ws1}) becomes
\begin{subequations}\label{split1}
\begin{align}
\pa_t v^\mathrm{c} & =  \lambda^\mathrm{c} v^\mathrm{c}
  + \pa_\vt p^\mathrm{c}_\mathrm{mf} 
\mathcal{N}(v^\mathrm{c},v^\mathrm{s}),\label{vceq} \\ 
\pa_t v^\mathrm{s} & = \Lambda^\mathrm{s} v^\mathrm{s}
 + P^\mathrm{s}_\mathrm{mf} 
\mathcal{N}(v^\mathrm{c},v^\mathrm{s}),\label{vseq}
\end{align}
\end{subequations}
where $\mathcal{N}$ is a smooth nonlinear map from
$H^{m+2}_\mathrm{ul}\times H^{m+2}_\mathrm{ul}$ into $H^m_\mathrm{ul}$
for every $m\geq1$.  
% and
% \[
% \mathcal{N}^\mathrm{c}(v^c,v^s):= \pa_\vt p^\mathrm{c}_\mathrm{mf} 
%         \mathcal{N}(v^c,v^s)
% \]
% maps $H^{m+2}_\mathrm{ul}\times H^{m+2}_\mathrm{ul}$ into 
% $H_\mathrm{ul}^m$ for each $m\geq0$. 

%%%%%%%%%%%%%%%%%%%%%%%%%%%%%%%%%%%%%%%%%%%%%%%%%%%%%%%%%%%%%%%%%%%%%%%%%%%%

\subsection{The moving frame}\label{mfr-sec}
To prove Theorems \ref{thstabil} and \ref{thmix1} 
we want to set up renormalization processes based on \reff{split1}.  
For this we need to remove the $\CO(\ell)$ terms in 
$\lam_1(\ell)=\ri(c_p-c_g)\ell-\al\ell^2+\CO(\ell^3)$. 
Therefore we define 
\begin{gather}
\bpm u^c\\u^s\epm=\CJ^{-1}\bpm v^c\\v^s\epm\quad\text{via}\quad 
\bpm \uti^c\\\uti^s\epm(\vt,\ell,t)=
\er^{\ri(c_p-c_g)\ell t}\bpm \vti^c\\ \vti^s\epm (\vt,\ell,t). 
\label{mframe}
\end{gather}
This yields 
\begin{subequations}\label{split2}
\begin{align}
\pa_t \uti^c =&\lamti_g(\ell)\uti^c
+(\pa_\vt{+}\ri\ell/k)\tilde{p}_{{\rm mf}}^c
\CNti(\uti^c,\uti^s), \label{uceq} \\ 
\pa_t \uti^s =&\Lamti_g^s(\ell)\uti^s
 + \tilde{P}^\mathrm{s}_\mathrm{mf}\CNti(\uti^c,\uti^s),\label{useq}
\end{align}
\end{subequations}
where $\lamti_g(\ell)=\lambda_1(\ell){-}\ri(c_p{-}c_g)\ell$ and 
$\Lamti_g^s(\ell)=\Lamti(\ell)-\ri(c_p-c_g)\ell$. 
The factors $\er^{\pm \ri c\ell t}$ 
drop out of the nonlinearities since as multipliers they commute with 
the mode filters and 
\begin{align*}
(\uti^{*2})(\ell)\er^{-\ri c\ell t}&=((\er^{\ri c\ell t}\vti)^{*2})(\ell)
\er^{-\ri c\ell t}
=\int_m \er^{\ri c(\ell-m) t}\vti(\ell-m)\er^{\ri cm t}\vti(m)\dd m\, 
\er^{-\ri c\ell t}=(\vti^{*2})(\ell), 
\end{align*}
and similar for higher power convolutions. 

In general, \reff{mframe} does not correspond to a simple transform 
in $\vt$-space. However, if $\uti$ has the special form 
$\uti(t,\ell,\vt)=\ati(\ell,t)g(\vt)$ then, cf.\,\reff{sft}, 
\begin{align}
v(\vt,t)&
=\int_{-k/2}^{k/2} \er^{\ri\ell(\vt/k+(c_p-c_g)t)}\ati(t,\ell)g(\vt)\dd\ell
=\al(\vt/k+(c_p-c_g)t,t)g(\vt)\notag\\
&=\bigl[\al(x-c_gt,t)+\pa_\vt\al(x-c_gt,t)\phi_\vt(\vt,t)
+\hot\bigr]g(\vt).
\label{btraf} 
\end{align}
Thus, \reff{mframe} will be responsible for recovering the group speed 
in Theorems \ref{thstabil} and \ref{thmix1}, 
which motivates the index $g$ in \reff{split2}. On the other hand, completely 
transforming \reff{split1} to a comoving frame would make the 
linear part 
spatially and temporally periodic, and the subsequent 
analysis would require Floquet theory in time and thus be 
more complicated. 

The key features of \reff{split2} are the following. By construction, 
\begin{gather}
\lam_g(\ell)=(\lambda_1(\ell){-}\ri(c_p{-}c_g)\ell)=
-\al\ell^2+\CO(\ell^3).\label{lamg}
\end{gather}
We have 
\begin{gather}
(\pa_\vt{+}\ri\ell/k)\tilde{p}_{{\rm mf}}^c\CNti(\uti^c+\uti^s)
=\etati(\ell)\CNti^c(\uti^c,\uti^s),
\end{gather}
where $|\etati(\ell)|=C\ell$ and $\CN^c$ maps $H^{m+2}(n)\times H^{m+2}(n)$ 
into $H^s(n)$ for all $s\in\N$. In particular, by the 
calculations from \cite{DSSS05}, 
\begin{gather}
\etati(\ell)\CNti^c(\uti^c,\uti^s)
=\beta\ri \ell(\uti^c)^{*2}+\hot, 
\label{kenel}
\end{gather}
with $\beta=- \frac{1}{2}\omnl''(k)$, and where the higher order terms $\hot$ are 
discussed later. The spectrum 
of $\Lamti^s_g$ is left of $\re z<-\sig_0$, hence $\uti^s$ is 
linearly exponentially damped. Thus, heuristically, if 
for now we ignore $\uti^s$ 
and $\hot$ in \reff{kenel}, then, as explained in \S\ref{buri-sec} 
we have the following situations: 
in Theorem \ref{thstabil} and in Theorem \ref{thmix1} case (i) 
(with $\omnl''=0$), corresponding 
to Proposition \ref{prop1} cases (i) and (ii), respectively, 
the whole nonlinearity is irrelevant and we 
obtain Gaussian diffusive behavior of $\uti^c$;  for 
 case (ii) of Theorem \ref{thmix1} ($\omnl''\neq 0$), corresponding 
to Proposition \ref{prop1} case (iii), 
the dynamics are governed by the Burgers equation for $\uti^c$. 

The (unavoidable) drawbacks of the coordinates \reff{split2} 
are their relatively complicated derivation, and 
that \reff{split2} is quasi-linear while the original system \reff{rd1} 
is semi-linear. 

%%%%%%%%%%%%%%%%%%%%%%%%%%%%%%%%%%%%%%%%%%%%%%%%%%%%%%%%%%%%%%%%%%%%%%%%%%%%

\section{The results in Bloch wave space}\label{res2-sec}

To prove Theorems \ref{thstabil} and \ref{thmix1}, in \S\ref{rg-sec} 
we set up renormalization processes for \reff{split2} 
in Bloch space. For this we need 
Bloch spaces with regularity and weights in $\ell$. Thus we first collect 
a number of definitions 
and basic properties. We recall that 
$H^{m_2}(m_1){=}\{u\in L^2(\R) : \| u \|_{H^{m_2}(m_1)}{<}\infty\}$ with 
$\| u\|_{H^{m_2}(m_1)} = \| u \rho^{m_1} \|_{H^{m_2}(\R)}$, where 
$ \rho(x) = (1+|x|^2)^{1/2}$, and that 
$\CF$ is an isomorphism between $H^{m_2}(m_1)$ and $\Hhat^{m_1}(m_2)$, 
where the notation $\Hhat^{m_1}(m_2)=H^{m_1}(m_2)$ is used to indicate 
functions that live in Fourier space. 

Similarly, for $L>0$ and $m_1,m_2,b\ge 0$ define 
\begin{gather*}
\CB_L^{m_1}(m_2,b):=\{\vti\in
H^{m_1}((-Lk/2,Lk/2),H^{m_2}_{{\rm per}}((0,2\pi))) \ : \
\|\vti\|_{\CB_L^{m_1}(m_2,b)}<\infty\},\notag\\
\|\vti\|_{\CB_L^{m_1}(m_2,b)}^2=\sum_{\alpha\leq m_1}
\sum_{\beta\leq m_2} 
\|(\pa_\ell^\alpha \pa_\vt^\beta\vti)\rho^{b}
\|^2_{L^2((-Lk/2,Lk/2),L^2(\CT_{2\pi}))}. 
\end{gather*}
Here $\rho=\rho(\ell)$, i.e., 
we introduce a weight in the Bloch wave number $\ell$, and 
 the subscript $L$ indicates that the Bloch wave number varies in 
$[-kL/2,kL/2]$. 
 For fixed $L>0$ 
the weight $\rho$ is irrelevant since, due to the bounded wave number domain,  
all norms $\|\cdot\|_{\CB_L^{m_1}(m_2,b_1)}$ and 
$\|\cdot\|_{\CB_L^{m_1}(m_2,b_2)}$ 
are equivalent, but the constants depend on $b_1,b_2$ and $L$. 
The purpose of the weights is to take advantage of the ``derivative 
structure'' of the nonlinearity in the equation for $\uti^c$, see 
\reff{kenel}, and Lemma \ref{dxlem} below. 

Let $\CB^{m_1}(m_2,b):=\CB_1^{m_1}(m_2,b)$. Then 
$\CJ$ is an isomorphism between $H^{m_2}(m_1)$ and $\CB^{m_1}(m_2,b)$, with 
arbitrary $b\geq 0$, see, e.g., 
\cite[Lemma 5.4]{gs98}. We define the scaling operators 
$$
\CR_{1/L}:\CB^{m_1}(m_2,b)\ra\CB_L^{m_1}(m_2,b),\quad 
[\CR_{1/L}\vti](\vt,\ell)=\vti(\vt,\ell/L).
$$
Only $\ell$ is rescaled, and $\vt$ is not, and similar to \reff{mframe} this 
does in general not correspond to a simple rescaling of 
$v$. However, note that $\uti^c=\uti^c(\ell,t)$ does not depend 
on $\vt$, i.e., for $\uti^c$ Bloch space is identified with Fourier space, 
and in this case we have 
\huga{\label{bresc}
\CJ^{-1}(\CR_{1/L}\uti)=\CF^{-1}(\CR_{1/L}\uti)=L\CR_Lu, 
}
i.e., concentration at $\ell=0$ in Bloch space corresponds to 
spreading in $\vt$. Finally, 
\begin{gather}
\|\CR_{1/L}\vti\|_{\CB_L^{m_1}(2,b)}\leq CL^{b+1/2}\|\vti\|_{\CB^{m_1}(2,b)}, 
\label{re1}
\end{gather}
and, for $\uti,\vti\in \CB_L^{m_1}(m_2,0)$ with 
$m_1,m_2\geq 1/2$ and $\ell\in(-L/2,L/2)$, 
\begin{align}
\CR_{1/L}(\CR_L\uti*\CR_L\vti)(\ell,x)
&=\int_{-1/2}^{1/2} \uti(\ell-Lm,x)v(Lm,x)\dd m
\notag\\
&=L^{-1}\int_{-L/2}^{L/2}\uti(\ell-m,x)\vti(m,x)\dd m
=:L^{-1}(\uti*_L\vti)(\ell,x). \label{bcone1}
\end{align}
This will be used to express the rescaled nonlinear terms, where 
henceforth we will drop the subscript $_L$ in $*_L$.

To recall the heuristics, as a model 
for Theorems \ref{thstabil} and \ref{thmix1}(i) (in which the 
nonlinearities are completely irrelevant),  
consider the Fourier transformed 
version of $\pa_t u=\al\pa_x^2 u$, $u_{t=1}=u_0$, i.e., 
$\pa_t \uti=-\al\ell^2 \uti$, 
which is solved by $\uti(\ell,t)=\er^{-(t-1)\al\ell^2}\uti(\ell,1)$. 
Then, for any $c\in\C$, or more specifically $c\in\R$ since we consider 
real valued functions $u$, $\fti_c(\ell)=c\er^{-\al\ell^2}$ is a fixed 
point of the renormalization map 
\begin{gather}
\CG^{(1)}:\uti\mapsto \er^{-\al\ell^2(1-1/L^2)}\CR_{1/L}\uti.
\label{rmap1}
\end{gather}
Moreover, for $L>1$ being sufficiently large, 
this line of fixed points is attractive in $H^2(2)$. 
To see this, write $\uti(\ell)=\fti_c(\ell)+\gti(\ell)$ with 
$\gti(0)=0$. Then, using $|\gti(\ell)|\leq (\ell/L)\|\pa_\ell\gti\|_{C^0_b}$ 
(by the mean value theorem) and $H^2\hookrightarrow C^1$ we obtain 
\begin{gather}
\|\er^{-\al\ell^2(1-1/L^2)}\CR_{1/L}\gti\|_{H^2(2)}^2
\leq CL^{-1}\|\gti\|_{H^2(2)}.\label{r1contr}
\end{gather}
Thus, $\vti(\ell,t)=\uti(t^{-1/2}\ell,t)\ra \fti_c(\ell)$ as $t\ra\infty$ 
is the expected scaling for $\uti^c$ in Theorem \ref{thmix1}(i). 
Theorem \ref{thmix1}(ii)  
is also based on \reff{rmap1} but we have a nonlinear correction 
to the asymptotic profile as explained in \S\ref{buri-sec}. 

Similarly, for any $c\in\R$, $\gti_c(\ell)=\ri c\ell\er^{-\al\ell^2}$ 
is a fixed point of the renormalization map 
\begin{gather}
\CG^{(2)}:\uti\mapsto \er^{-\al\ell^2(1-1/L^2)}L\CR_{1/L}\uti, 
\label{rmap2}
\end{gather}
and again this line of fixed points is attractive in $H^3(2)\cap X_0$, 
where $X_0$ consists of functions with zero mean. 
For this write $\uti(\ell)=\gti_c(\ell)+\hti(\ell)$ with 
$\pa_\ell\hti(0)=0$ and use $|\hti(\ell/L)|\leq (\ell/L)^2
\|\pa_\ell^2\hti\|_{C^0_b}$ and $H^3\hookrightarrow C^2$. Thus, 
 $\vti(\ell,t)=t^{1/2}\uti(t^{-1/2}\ell,t)\ra \gti_c(\ell)$ as $t\ra\infty$ 
is the expected scaling for $\uti^c$ in Theorem \ref{thstabil}. 
The need for $\uti\in C^2$ also explains 
the higher weight in $x$ in Theorem \ref{thstabil}.

\begin{Theorem} \label{thstabil-b}{\bf [Diffusive stability].}
Let $ u_0(\cdot;k)$ be a spectrally 
stable wave--train and $ b \in (0,1/2) $. There exist $\eps,C>0$ such that 
if $\|(\uti^c,\uti^s)|_{t=1}\|_{H^3(2)\times\CB^3(2,2)}\leq \eps_1$ 
and $\uti^c(0,1)=\frac 1 {2\pi k}\int u^c(\vt,1)\dd \vt=0$, 
then the solution $(\uti^c,\uti^s)$ to \reff{split2} exists for all 
$t\geq 1$, and there exists a $\psilimti\in\R$ such that 
\begin{align}
\|t^{1/2}\uti^c(t^{-1/2}\ell,t)
-\ri\psilimti\ell\er^{-\al\ell^2}\|_{H^3(2)}&\leq 
C t^{-1/2+b}, \label{uc}\\
\|t^{1/2}\uti^s(t^{-1/2}\ell,t)\|_{\CB_{\sqrt{t}}^3(2,2)}&
\leq C t^{-1/2+b}. \label{us} 
\end{align}
\end{Theorem}

\begin{Theorem} \label{thmix1-b}{\bf [Diffusive mixing of phases].}
Let $ u_0(\cdot;k)$ be a spectrally 
stable wave--train and $ b \in (0,1/2) $. 
There exist $\eps,C>0$ such that for $|\phi_d|\leq\eps$ 
the following holds. 

(i) Assume that  $\beta=- \frac{1}{2}\omnl''(k)=0$, 
$\|\uti^c(\ell,1)\|_{H^2(2)}\leq \eps$ with $\uti^c(0,1)=\phi_d/(2\pi k)$, 
$\|\uti^s(\cdot,1)\|_{\CB^2(2,2)}\leq \eps$ and 
$\tilde{P}^s\uti^s(\cdot,1)=\uti^s(\cdot,1)$. 
Then the solution $(\uti^c,\uti^s)$ to \reff{split2} exists for all 
$t\geq 1$, and 
\begin{align} 
\|\CR_{1/\sqrt{t}}\uti^c(\ell,t)-\uti^c_*(\ell)\|_{H^2(2)}&\leq 
Ct^{-1/2+b}, \label{n27b}\\
\|\CR_{1/\sqrt{t}}\uti^s(\ell,t)\|_{B_{\sqrt{t}}^2(2,2)}&\leq Ct^{-1/2+b}\ , 
\label{n27b2} 
\end{align}
where 
$\uti_*^c(\ell)=\phi_d\er^{-\al\ell^2}$.  

(ii) If $\beta=- \frac{1}{2}\om''(k)\neq 0$ then the same result holds with 
$\uti_*^c(\ell)$ replaced by 
\begin{gather} \label{n27c}
\uti_*^c(\ell)=
\CF\left(\frac{\sqrt{\al}}\beta\frac{z\er^{-\vt^2/(k^2\al)}}
{1+z\erf(\vt/\sqrt{k\al})}\right)(\ell), 
\end{gather}
where $\ln(1+z)=\frac{\beta}{\al}\phi_d$. 
\end{Theorem}

\noindent
%Note that $\uti^c(0,1)=\phi_d/(2\pi k)$ and $\tilde{P}^s\uti^s(\cdot,1)=\uti^s$ 
%imply that and $\phi_\pm\in [0,2\pi)$ are arbitrary but the condition 
%$t_1>t_0$ for some possibly large $t_0$ ensures that the initial 
%phase interface is sufficiently flat. 
Before proving Theorems \ref{thstabil-b} and \ref{thmix1-b} we 
show that they imply Theorems \ref{thstabil} and \ref{thmix1}.

\noindent
{\bf Proof of Theorem \ref{thstabil}.} \ \ 
Given initial data in the form \reff{icdecom1}
from Theorem \ref{thstabil}, i.e., 
\begin{gather*}
u(x,t)|_{t=0}=u_0(\th-\th_0+\phi_0(x);k)+v_0(x)\quad\text{with}\quad 
\|\phi_0\|_{H^3(3)}, \| v_0 \|_{H^3(2)} \le \eps, 
\label{icdecom1b}
\end{gather*} 
we first need to extract $(\uti^c,\uti^s)|_{t=1}$ 
and show that they fulfill the assumptions of Theorem \ref{thstabil-b}. 
Then we translate \reff{uc},\reff{us} back into $(\phi,v)$ coordinates. 

Thus, as explained in Remark \ref{inter-rem}, let 
$u(x,t)|_{t=1}=u_0(\vt;k(1+\pa_\vt\phi_0))+w_0(\vt)$ 
with $\th=\vt-\phi_0(\vt)$ 
%where $\vt-\phi_0(\vt)$ is the inverse of $\th\mapsto\th+\check{\phi}(k\th)$, 
and 
$$
w_0(\vt):=u_0(\vt;k)-u_0(\vt;k(1+\pa_\vt\phi_0))+v_0(x). 
$$
W.l.o.g.\ assume that $(1-P^{{\rm s}})w_0=0$, otherwise redefine 
$\phi_0=p_{{\rm mf}}^{{\rm c}}\phi_0$. This fixes the 
non-uniqueness in \reff{icdecom1}. 
Also, $\phi_0\in H^m(3)$ for all $m\in\N$ due to the compact support 
of $\tilde{\phi}_0$, %, hence $w_0(\vt)$ is small in $H^3(2)$.  
and, with $\psi_0=\pa_\vt\phi_0$, $J$ from \reff{jdef} and $\CS$ 
from \reff{sdef}, 
\begin{gather}
(\uti^c,\uti^s)|_{t=1}=J\CS^{-1}(\psi_0,w_0)=J(\psi_0,w_0-S_1\psi_0)
\label{ap1}
\end{gather}
is well defined and fulfills 
$\|(\uti^c,\uti^s)\|_{H^3(2)\times\CB^3(2,2)}\leq C_1 \eps$ and 
$\uti^c(0,1)=\frac 1 {2\pi k}\int u^c(\vt,1)\dd \vt=0$. 

We now use \reff{uc},\reff{us} to recover Theorem \ref{thstabil}. 
Using \reff{sft2} and 
$\ds 
{\cal F}_1^{-1}(\ri\ell\er^{-\al\ell^2})(\vt)=-\frac 1{\sqrt{4\pi\al}}
\frac{\vt}{2\al}\er^{-\vt^2/(4\al)}
$
we have 
\begin{align}
t\CR_{t^{1/2}}&u^c(\vt,t)-\psilim\vt\er^{-\vt^2/(4k^2\al)}
=\CF^{-1}\left[t^{1/2}\CR_{t^{-1/2}}\uti_c(\ell,t)
-\ri\psilimti\ell\er^{-\al \ell^2}\right](\vt), \notag
\end{align}
where 
$\ds \psilim=-\frac{\psilimti}{\sqrt{4\pi \al}}\frac{1}{2\al k}$, and from 
$c_1\|\uhat\|_{H^n(m)}\leq \|u\|_{H^m(n)}\leq c_2\|\uhat\|_{H^n(m)}$ 
we obtain 
$$
\|t\CR_{t^{1/2}}u^c(\vt,t)-\psilim\vt\er^{-\vt^2/(4k^2\al)}\|_{H^2(3)}
\leq Ct^{-1/2+b}. 
$$
Then, with 
$$
\psi(\vt,t)=u^c(\vt+k(c_p-c_g)t,t)\quad\text{and}\quad 
w(\vt,t)=\CS_1u^c(\vt/k+(c_p-c_g)t,t)+u^s(\vt,t)
$$
we obtain, in $L^\infty$, 
\begin{align}
\phi(\vt,t):=&\int_{-\infty}^\vt \psi(\xi,t)\dd \xi
=-\frac{2k^2\psilim\al}{\sqrt{t}}\int_{-\infty}^{\vt+k(c_p-c_g)t}
\left(-\frac{\xi}{2\al k^2 t}\er^{-\xi^2/(4k^2\al t)}\right) 
\dd \xi+\CO(t^{-1})\notag\\
=&-2t^{-1/2}k^2\psilim\al\er^{-(\vt+k(c_p-c_g)t)^2/(4\al k^2 t)}+\CO(t^{-1}), 
\label{phiinf}
\end{align}
i.e., $\philim=-4k^2\sqrt{\al^3\pi}\psilim$. Also $w(\vt,t)=\CO(t^{-1})$
since 
$\supp\tilde{\CS}_1\subset\{\ell_1/8<|\ell|<\ell_1/4\}$ and 
\begin{align}
\|u^s(\vt,t)\|_{L^\infty}&=\|J^{-1}\uti^s(\cdot,\cdot,t)(\vt)\|_{L^\infty}
=\left\|\int_{-k/2}^{k/2}\er^{\ri\ell\vt/k}\uti^s(\vt,\ell,t)
\dd \ell\right\|_{L^\infty}\notag\\
&=\left\|t^{-1/2}\int_{-kt^{1/2}/2}^{kt^{1/2}/2}\er^{\ri\ell t^{-1/2}\vt/k}
\uti^s(\vt,t^{-1/2}\ell,t)(1+\ell^2)^2(1+\ell^2)^{-2}
\dd \ell\right\|_{L^\infty}\notag\\
&\leq Ct^{-1/2}\|\CR_{t^{-1/2}}\uti^s(\cdot,\cdot,t)\|_{B_{\sqrt{t}}^3(2,2)}
\leq Ct^{-1+b}.\label{usinf}
\end{align}
Finally, 
\begin{align}
\vt=&\th-\th_0+\phi(\vt,t)=\th-\th_0-2k^2t^{-1/2}\psilim\al
\er^{-(\vt+k(c_p-c_g)t)^2/(4\al k^2 t)}+\CO(t^{-1})\notag\\
=&\th-\th_0-2k^2t^{-1/2}\psilim\al
\er^{-(x-c_gt)^2/(4\al t)}+\CO(t^{-1})\notag
\end{align}
using $\th=kx-\om t=k(x-c_p t)$ and the implicit function theorem.\qed 

\vs{2mm}
\noindent
{\bf Proof of Theorem \ref{thmix1}.} \ \
First, assume that $\beta=0$. As above we write 
$u(x,t)|_{t=0}=u_0(\vt;k(1+\pa_\vt\phi_0)+w_0(\vt)$ 
with $w_0(\vt):=u_0(\vt;k)-u_0(\vt;k(1+\pa_\vt\phi_0))+v_0(x)$, 
and where now 
$\|\phi_0'(\cdot)\|_{H^2(2)}\leq \eps$ and $\phi_0(\vt)\ra\phi_\pm$ as 
$\vt\ra\pm\infty$. Again, w.l.o.g.~assume that $(1-P^{{\rm s}})w_0=0$. 
Then, for 
$$
\psi_0(\vt)=\pa_\vt\phi_0(\vt)=u^c(\vt,1) 
$$ 
we obtain $2\pi k\uti^c(0,1)=\int u^c(\vt,1)\dd\vt=\phi_d$, and 
Theorem \ref{thmix1-b} applies to 
$$
(\uti^c,\uti^s)|_{t=1}=J\CS^{-1}(\psi_0,w_0)=J(\psi_0,w_0-S_1\psi_0).
$$ 
Thus, with 
$\psi(\vt,t)=u^c(\vt+k(c_p-c_g)t,t)$ and 
$w(\vt,t)=\CS_1u^c(\vt+k(c_p-c_g)t,t)+u^s(\vt,t)$ we obtain, in $L^\infty$, 
\begin{align*}
\phi(\vt,t):=&\phi_-+\int_{-\infty}^\vt \psi(\xi,t)\dd \xi
=\phi_-+(\phi_+-\phi_-)
\frac{1}{\sqrt{4\pi k^2t}}\int_{-\infty}^{\vt+k(c_p-c_g)t}
\er^{-\xi^2/(4k^2\al t)}
\dd \xi+\CO(t^{-1})\\
=&\phi_-+(\phi_+-\phi_-)\erf((x-c_gt)/\sqrt{4\al t})+\CO(t^{-1})
\end{align*}
and $w(\vt,t)=\CO(t^{-1/2+b})$ as in \reff{phiinf} and \reff{usinf} above. 
Hence 
$$
\vt=\th+\phi(\vt,t)=\th+\phi_-+(\phi_+-\phi_-)
\erf(x-c_gt)/\sqrt{4\al t})+\CO(t^{-1/2})
$$
and by shifting the $\CO(t^{-1/2})$-part to $v$ we obtain part (i) in Theorem 
\ref{thmix1}. Part (ii) with $\beta\neq 0$ works in the same way. \qed

%%%%%%%%%%%%%%%%%%%%%%%%%%%%%%%%%%%%%%%%%%%%%%%%%%%%%%%%%%%%%%%%%%%%%%%%%%%%

\section{Renormalization}\label{rg-sec} 

We first prove Theorem \ref{thstabil-b}; 
the minor modifications needed to prove Theorem \ref{thmix1-b}(i) 
are then explained in \S\ref{p2-sec}, while the changes for the 
slightly more complicated proof of 
Theorem \ref{thmix1-b}(ii) are explained in \S\ref{p3-sec}. 

%%%%%%%%%%%%%%%%%%%%%%%%%%%%%%%%%%%%%%%%%%%%%%%%%%%%%%%%%%%%%%%%%%%%%%%%%%%%

\subsection{The rescaled systems}\label{re-sec}
Based on the asserted behavior 
$t^{1/2}\uti^c(t^{-1/2}\ell,t)\ra \ri\psilim\ell \er^{-\al\ell^2}$ we 
introduce, for $n\in\N$ and $L>1$ chosen sufficiently large below, 
the variables 
\begin{align} 
\uti_n^c(\kap,\tau)&:=
L^n\uti^c(\kap/L^n,L^{2n}\tau)= 
L^n[\CR_{L^{-n}}\uti^c](\kap,L^{2n}\tau),\\
\uti_n^s(\vt,\kap,\tau)&:=L^n\uti^s(\vt,\kap/L^n,L^{2n}\tau)
=L^n[\CR_{L^{-n}}\uti^s](\vt,\kap,L^{2n}\tau). 
\end{align}
Then $(\uti_n^c,\uti_n^s)$ fulfill
\begin{subequations}\label{splitn}
\begin{align}
\pa_\tau \uti_n^c(\kap,\tau)-\lamti_{g,n}(\kap)\uti_n^c(\kap,\tau)
&=L^{3n}\CNti_n^c(\uti_n^c,\uti_n^s)( \kap,\tau),\label{unce}\\
\pa_\tau \uti_n^s(\vt,\kap,\tau)-
\Lamti_{g,n}\uti_n^s(\vt,\kap,\tau)
&=L^{3n}\CNti_n^s(\uti_n^c,\uti_n^s)\vt,\kap,\tau),\label{unse}
\end{align}
\end{subequations}
where 
\begin{align*}
&\lamti_{g,n}(\kap)=L^{2n}\lamti_g(\kap/L^{n}),\quad 
\Lamti_{g,n}=L^{2n}\CR_{L^{-n}}\Lamti_{g}\CR_{L^n},\\ 
&\CNti_n^c(\uti_n^c,\uti_n^s)=
\etati(\kap/L^{n})\CR_{L^{-n}}\CNti^c(L^{-n}\CR_{L^n}\uti_n^c,L^{-n}\CR_{L^n}
\uti_n^s),\\
&\CNti_n^s(\uti_n^c,\uti_n^s)=
\CR_{L^{-n}}\pmfsti\CNti(L^{-n}\CR_{L^n}\uti_n^c,L^{-n}\CR_{L^n}\uti_n^s). 
\end{align*}
Except for the different scaling due to $\uti_n^c(0,\ell)=0$, 
\reff{splitn} has a very similar structure as, e.g., 
\cite[eq.(30)]{gs96a} or \cite[eq.(3.2)]{ue04ibl}. 
Thus, similar to \reff{r0}, we shall consider the following iteration: 
\begin{gather} 
\begin{split}
\text{solve \reff{splitn} for } \tau\in I:=[L^{-2},1]\text{ with initial data }
\bpm \uti_n^c\\\uti_n^s\epm
(\vt,\kap,L^{-2})=L\bpm\uti_{n-1}^c\\ \uti_{n-1}^s\epm
(\vt,\kap/L,1).
\end{split}
\label{r1}
\end{gather}
%As phase space for \reff{splitn} we choose $\CX_n\times\CX_n$ with 
%$\CX_n=\CB_{L^n}^3(2,2)$, where for 
%$\uti_n^c$ we can identify $\CX_n$ with the 
%Fourier space $H^3(2)$. Moreover, 
%$\supp\uti_n^c\subset\{|\kap|\leq L^n\ell_1\}$. 

Formally, \reff{splitn} is solved by the variation of constant formula, i.e. 
\begin{subequations}\label{vdkn}
\begin{align}
\uti_n^c(\kap,\tau)=&\ \er^{(\tau-1/L^2)\lamti_{g,n}(\kap)}
\uti_n^c(\kap,1/L^2)\notag\\
&+\int_{1/L^2}^\tau \er^{(\tau-s)\lamti_{g,n}(\kap)}
L^{3n}\CNti_n^c(\uti_n^c,\uti_n^s)( \kap,s)\dd s,\label{uncsol}\\
\uti_n^s(\vt,\kap,\tau)=&\ \er^{(\tau-1/L^2)\Lamti_g\CR_{L^n}}
\uti_n^s(\vt,\kap,\tau)\notag\\
&+\int_{1/L^2}^\tau\er^{(\tau-s)\Lamti_{g,n}}
L^{3n}\CNti_n^s(\uti_n^c,\uti_n^s)(\vt,\kap,s)\dd s.\label{unssol}
\end{align}
\end{subequations}
However, \reff{vdkn} can not be used to construct the solution since 
\reff{splitn} is a quasi-linear system, as it can be seen from 
$\CN:H^{m_2}(m_1)\times H^{m_2}(m_1)\ra H^{m_2-2}(m_1)$ in \reff{split1}. 
To solve \reff{splitn} we use maximal 
regularity methods \cite{lm68} for parabolic equations 
in (weighted) Sobolev spaces as in \cite{ue04nsb}. 
A posteriori, \reff{vdkn} can then be used to estimate the solutions. 
Thus, we first note some properties of the linear semigroups and the 
nonlinearities in \reff{vdkn}, and then explain how to obtain 
local existence for \reff{splitn}. 

\subsection{Estimates on the linear semigroups and 
the nonlinearities} 
\label{nlest-sec}
We shall need some detailed estimates on the linear semigroups and 
the nonlinear terms in \reff{vdkn}. The idea is to exploit the 
derivative-like structure 
in the Bloch wave number $\kap$ of $\CNti^c$ as expressed in \reff{kenel} 
by relaxing the weight, and to regain the weight using 
$\er^{(\tau-\tau')\lamti_{g,n}}$. 
Thus, from this point on, the weights in $\kap$ become important. 

\begin{Lemma}\label{lsg-lem}
There exists a $C>0$ such that for all $L>1$ we have 
\begin{align}
\|\er^{(\tau-\tau')\lamti_{g,n}}\uti_n^c\|_{B_{L^n}^3(2,2)}
\leq &C\max\{1,(\tau-\tau')^{-b/2}\}\|\uti_n^c\|_{B_{L^n}^3(2,2-b)},
\label{csge}\\
\|\er^{(\tau-\tau')\Lamti_{g,n}}\uti_n^s\|_{B_{L^n}^3(2,2)}
\leq &C\max\{1,(\tau-\tau')^{-m_2/2}\}\er^{-\ga_0L^{2n}(\tau-\tau')}
\|\uti_n^s\|_{B_{L^n}^3(2-m_2,2)}.
\label{ssge}
\end{align}
\end{Lemma}
{\bf Proof.} Equation \reff{csge} holds since the real part of 
$\lamti_{g,n}(\kap)=L^{2n}\lamti_g(\kap/L^n)
{=}-\al\kap^2+\CO(\kap^3)$ is bounded from above by the 
parabola $-\al_0\kap^2$, 
while \reff{ssge} holds since 
$\Lamti_{g,n}$ is a relatively bounded perturbation of 
$L^{2n}(\pa_\vt+\ri\kap/L^n)^2$ and by construction has spectrum left 
of $-L^{2n}\gamma_0$.\qed

The following lemma transfers the fact that derivatives give higher powers 
of $L^{-1}$ upon rescaling to general convolution operators 
with a ``derivative--like'' structure. 
\begin{Lemma}\label{dxlem}
Let $m_1\in\N$, $\ga\geq 0$, and 
$\tilde{K}\in C^{m_1}_{b}([-1/2,1/2)^2,H^{2}(\CT_{2\pi}))$ with 
$
\|\tilde{K}(\kap{-}\ell,\ell)\|_{H^{2}(\CT_{2\pi})}
\leq C(|\kap{-}\ell|{+}|\ell|)^\gamma. 
$
Then 
$$ 
(\vti,\wti)\mapsto (\CM_{1/L} K)(\vti,\wti)(\kap,\vt):=
\int_{-L/2}^{L/2} \bigl[\CR_{1/L}\tilde{K}\bigr](\kap-\ell,
\ell,\vt)\vti(\kap,\vt)\wti(\kap-\ell,\vt)\dd\ell
$$
defines a bilinear mapping 
$(\CM_{1/L} K):\CB_L^{m_1}(2,2)\times \CB_L^{m_1}(2,2)\ra 
\CB_L^{m_1}(2,2)$, and there exists a 
$C>0$ such that for all $L>1$ we have
$$
||(\CM_{1/L} K)(\vti,\wti)||_{\CB_L^{m_1}(2,2-\gamma)}
\leq CL^{-\min\{\gamma,1\}}||\vti||_{\CB_L^{m_1}(2,2)} 
||\wti||_{\CB_L^{m_1}(2,2)}.
$$
\end{Lemma}
{\bf Proof.} This holds due to $\sup_\ell\left|\frac 
{\ell^\ga L^{-\ga}}{(1+\ell^2)^{\ga/2}}\right |\leq CL^{-\ga}$. \qed

\begin{Lemma}\label{non-lem} 
Let $\|(\uti_n^c,\uti_n^s)\|_{[B^3_{L^n}(2,2)]^2}\leq R_n\leq 1$. 
There exists a $C>0$ such that 
\begin{align}
L^{3n}\|\CNti_n^c(\uti_n^c,\uti_n^s)\|_{B^3(2,1)}&\leq CL^{-n}R_n^2. 
\label{ntice}
\end{align}
The term 
$L^{3n}\CNti_n^s$ can be split according to the number of $\vt$ derivatives 
as $L^{3n}\CNti_n^s=\CNti_{n,0}^s+\CNti_{n,1}^s+\CNti_{n,2}^s$ such that 
\begin{align}
\|\CNti_{n,i}^s\|_{B^3(2-i,2)}&\leq CR_n^2.
\label{ntise}
\end{align}
\end{Lemma}
{\bf Proof.} 
We write $L^{3n}\CNti_n^c=s_1+s_2$, where, as 
explained in \S\ref{mfr-sec}, the lowest order term $s_1$ in 
$L^{3n}\CNti^c_n(\uti_n^c,\uti_n^s)$ reads 
\begin{gather} 
s_1(\kap)=L^{3n}\ri\beta\frac{\kap}{L^n}
\CR_{L^{-n}}(L^{-n}\CR_{L^n}\uti_n^c)^{\star 2}(\kap),
\label{s1def}
\end{gather}
cf.\,\reff{kenel}. This yields $\|s_1\|_{\CB^3(2,1)}\leq CL^{-n}R_n^2$ by 
direct calculation. The remaining terms $s_2$ can be estimated in a 
similar way using Lemma \ref{dxlem} and taking into account the 
finite support of 
$(\pa_\vt{+}\ri\ell/k)\tilde{p}_{{\rm mf}}^c\CNti(\uti^c+\uti^s)$ 
in Fourier space. 

This does not work for $L^{3n}\CNti_n^s(\uti_n^c,\uti_n^s)$. However, 
here we do not need an additional factor $L^{-n}$, and \reff{ntise} 
simply follows by checking the number of derivatives in $\CN$ and 
using \reff{bcone1}.  \qed

%%%%%%%%%%%%%%%%%%%%%%%%%%%%%%%%%%%%%%%%%%%%%%%%%%%%%%%%%%%%%%%%%%%%%%%%%%%%

\subsection{Local existence}\label{loc-ex-sec}
Since \reff{split1} and hence \reff{splitn} is quasilinear we cannot 
combine Lemma \ref{lsg-lem} 
and Lemma \ref{non-lem} to directly show local existence for \reff{splitn} 
via \reff{vdkn}. Instead we use maximal regularity theory from \cite{lm68}. 
For $I=(\tau_0,\tau_1)$ and $r,s\geq 0$ let 
$$
H^{r,s}(I,m_1)=L^2(I,H^r(m_1))\cap H^s(I,L^2(m_1). 
$$
Since \reff{split1} is a parabolic problem these spaces only occur 
with $s=r/2$ and we set $K^{m_2}(I,m_1)=H^{m_2,m_2/2}(m_1)$. 
Then, for any given weight $b>0$, Bloch transform is an isomorphism 
between $K^{m_2}(I,m_1)$ and 
$$
\Kti^{m_1}(I,m_2,b)=L^2(I,B^{m_1}(m_2,b))\cap H^{m_2/2}(I,B^{m_1}(0,b)). 
$$
Similarly, for every $n$, let 
$$
\Kti_{L^n}^{m_1}(I,m_2,b):=\CR_{1/L^n}\Kti^{m_1}(I,m_2,b):=
L^2(I,B_{L^n}^{m_1}(m_2,b))\cap H^{m_2/2}(I,B_{L^{n}}^{m_1}(0,b)), 
$$
i.e., the subscript $L^n$ again indicates that the Bloch wave number 
varies in $[-kL^n/2,kL^n/2)$. From \reff{re1} we have 
\huga{\label{re17}
\|\CR_{1/L^n}\uti\|_{\Kti_{L^n}^{m_1}(I,m_2,b)}
\le CL^{n(b+1/2)}\|\uti\|_{\Kti^{m_1}(I,m_2,b)}
}
Recall that for each $n$ the weight $b$ in $\kap$ 
gives an equivalent norm in $\Kti_{L^n}^{m_1}(I,m_2,b)$, but 
the constants depend on $n$.  
%The idea (cf.~Lemma \reff{non-lem}) 
%is to trade some weight in $\kap$ for factors of $L^{-\ga}$ 
%in the diffusive part, and 
%to recover the weight from the linear semigroup. 
We also need subspaces of functions that vanish sufficiently fast at $\tau_0$, 
and define  
\begin{gather*}
_0K^{m_2}(I,m_1):=\{v\in K^{m_2}(I,m_1): 
\pa_\tau^j v(\cdot,\tau_0)=0\text{ for } 
j\in\N, j<m_2/2-1/2\},\\
_0\Kti^{m_1}_{L^n}(I,m_2,b):=\{v\in \Kti^{m_1}_{L^n}(I,m_2,b): 
\pa_\tau^j \vti(\cdot,\cdot,\tau_0)=0\text{ for } 
j\in\N, j<m_2/2-1/2\}
\end{gather*}

We set 
$$I=(L^{-2},1),
$$ and for 
$(\uti_n^c,\uti_n^s)|_{\tau=L^{-2}}\in [B_{L^n}^3(2,2)]^2$ 
construct solutions $(\uti_n^c,\uti_n^s)\in [\Kti_{L^n}^3(I,3,2)]^2$ 
to \reff{splitn}. 
Note again that for $\uti_n^c$ we can identify Bloch space with 
Fourier space such that in fact $\uti_n^c\in K^3(I,2)$ (in the 
Fourier sense) 
with $\supp\,\uti^c_n(\tau)\subset I_n=\{|\kap|\leq L^n\ell_1/4\}$. 
We abbreviate \reff{splitn} as $\CL_n\Uti_n=\CNti_n(\Uti_n)$, where 
\begin{gather}
\CL_n\Uti_n=\bpm 
\pa_\tau \uti_n^c(\kap,\tau)-\lamti_{g,n}(\kap)\uti_n^c(\kap,\tau)\\
\pa_\tau \uti_n^s(\vt,\kap,\tau)-
\Lamti_{g,n}\uti_n^s(\vt,\kap,\tau), 
\epm, 
\end{gather}
and, for $m_2\geq 2$, we first consider the linear inhomogeneous version 
of \reff{splitn} with zero initial data, i.e., 
\begin{gather}
\CL_n \Uti_n(\tau)=\CNti_n(\tau),\quad 
\CNti_n\in{}[_0\Kti_{L^n}^3(I,m_2-2,2)]^2, \quad \Uti_n|_{\tau=L^{-2}}=0, 
\label{splitni}
\end{gather}
where moreover for the first 
component $\CNti_n^c$ of $\CNti_n=(\CNti_n^c,\CNti_n^s)$ we assume 
\huga{\label{nnca}\CNti_n^c\in K^3(I,2) \text{ (in the 
Fourier sense), and } 
\supp\,\CNti^c_n(\tau)\subset I_n=\{|\kap|\leq L^n\ell_1/4\}.
}

\begin{Lemma}\label{lsg1-lem} There exists a $C>0$, independent of $n\in\N$,  
such that for all $\CNti_n\in{}[_0\Kti_{L^n}^3(I,m_2-2,2)]^2$ which 
fulfill \reff{nnca}
there exists a unique solution of \reff{splitni} with 
\begin{gather}\label{line23}
\|\Uti_n\|_{\Kti_{L^n}^3(I,m_2,2)}\leq C 
\|\CNti_n\|_{_0\Kti_{L^n}^3(I,m_2-2,2)}.
\end{gather}
\end{Lemma}
{\bf Proof.} %Writing $\CNti_n=(\CNti_n^c,\CNti_n^s)$, 
The first component 
$\pa_\tau \uti_n^c(\kap,\tau)-\lamti_{g,n}(\kap)\uti_n^c(\kap,\tau)=
\CNti_n^c$ is independent of $\vt$ and thus can be solved by the 
variation of constant formula using \reff{csge} (with $b=0$). For the second 
component we use resolvent estimates for the solution of 
$$
(\lam-\Lamti_{g,n})\uti_n^s=\CNti_n^s.
$$
There exists a $C>0$ such that for $m_2\ge 2$, 
$\CNti_n^s\in B^{m_1}_{L^n}(m_2-2,b)$ all $\lam\in\C$ with 
$\re\lam\ge 0$ we have 
\huga{\label{line24}
\|\uti_n^s\|_{B^{m_1}_{L^n}(m_2,b)}+|\lam|^{m_2/2}\|u\|_{B^{m_1}_{L^n}(0,b)}
\le C\biggl(\|\CNti_n^s\|_{B^{m_1}_{L^n}(m_2-2,b)}+|\lam|^{(m_2-2)}
\|\CNti_n^s\|_{B^{m_1}_{L^n}(0,b)}\biggr)
} 
Similar to Lemma \ref{lsg-lem} 
this holds since 
$\Lamti_{g,n}$ is a relatively bounded perturbation of 
$L^{2n}(\pa_\vt+\ri\kap/L^n)^2$ and by construction has spectrum 
left of $-L^{2n}\ga_0$. See also \cite[Appendix A.2]{ue04nsb} for an 
explanation of how to obtain resolvent estimates in 
weighted spaces. From \reff{line24} we obtain \reff{line23} 
by continuation of $\CNti_n^s$ for $\tau\in\R$, 
Laplace transform, and the Paley--Wiener Theorem. In fact, 
in \reff{line24} we could choose $\lam$ to the right 
of $-L^{2n}\ga_0$, but  $\re\lam\ge 0$ is enough 
to show \reff{line23} with $C$ independent of $n$. \qed

We denote the solution operator of \reff{splitni} by $_0\CL_{n}^{-1}$. 
To solve the nonlinear problem we write $\Uti_n=\Vti_n+\Wti_n$ where 
$\Vti_n\in \Kti_{L^n}^3(\R,3,2)$ is a continuation of
$\Uti_n|_{\tau=L^{-2}}$, which exists due to \cite[Thm 4.2.3]{lm68}. 
Then $\Wti_n$ fulfills 
\begin{gather}
\CL_n\Wti_n=G_n(\Wti_n),\quad \Wti_n|_{\tau=L^{-2}}=0, 
\quad\text{where}\quad G_n(\Wti_n)=\CNti_n(\Vti_n+\Wti_n)-\CL_n\Vti_n.
\label{con0}
\end{gather}
The idea is to show that for $\Wti_n\in\,_0\Kti^3_{L^n}(I,3,2)$ we have 
$\Gti(\Wti_n)\in\,_0\Kti^3_{L^n}(I,1,2)$ and use Lemma \ref{lsg1-lem} and 
estimates on the nonlinearity to apply the contraction mapping 
theorem to 
\begin{gather}
\Phi(\Wti):=\, _0\CL_{n}^{-1}G_n(\Wti_n). 
\label{con1}
\end{gather}
We set 
\begin{gather}
\rho_n:=\|(\uti_{n}^c,\uti_{n}^s)|_{\tau=1}\|_{[B_{L^n}^3(2,2)]^2}
\label{rndef}
\end{gather}
and obtain the following local existence result, taking into 
account that $(\uti_n^c,\uti_n^s)(1/L^2)$ and 
$(\uti_{n-1}^c,\uti_{n-1}^s)$ are related by 
$(\uti_n^c,\uti_n^s)|_{\tau=L^{-2}}=
L\CR_{1/L}(\uti_{n-1}^c,\uti_{n-1}^s)|_{\tau=1}$ and hence, by \reff{re1}, 
$$
\|(\uti_n^c,\uti_n^s)|_{\tau=L^{-2}}\|_{B^3_{L^n}(2,2)}\leq 
CL^{7/2}\rho_{n-1}. 
$$
\begin{Lemma}\label{loc-ex1} There exist $C_1,C_2>0$, independent of $n$, 
such that the following holds. If $\rho_{n-1}\leq C_1L^{-7/2}$, 
then there exists a unique solution 
$(\uti_n^c,\uti_n^s)\in [\Kti_{L^n}^3(I,3,2)]^2$ to \reff{splitn} with 
\begin{gather}
\|(\uti_n^c,\uti_n^s)\|_{[\Kti_{L^n}^3(I,3,2)]^2}\leq C_2L^{7/2}\rho_{n-1}. 
\label{ee1}
\end{gather}
Moreover, for all $\tau_1>L^{-2}$ and any $m_2\in \N$ 
there exists a $C_3$, independent of $n$,  such that 
\begin{gather}
\|(\uti_n^c,\uti_n^s)\|_{[\Kti_{^n}^3((\tau_1,1),m_2,2)]^2}
\leq C_3L^{7/2}\rho_{n-1}.
\label{hreg1}
\end{gather}
\end{Lemma}
{\bf Proof.}  From standard Sobolev embeddings we have that $\CNti_n$ 
is a smooth mapping from $\Kti^3_{L^n}(I,3,2)$ to $\Kti^3_{L^n}(I,1,2)$, see 
also Lemma \ref{non-lem}.  
To show that 
$G_n(\Wti_n)$ in \reff{con0} is in $_0\Kti^3_{L^n}(I.1,2)$ we have to fulfill 
one compatibility condition, namely $G_n(\Wti_n)|_{\tau=L^{-2}}=0$, 
which holds by construction. For sufficiently small $\rho_{n-1}$, 
$\Phi$ is a contraction since $\CNti_n$ is quadratic and higher order. 
In particular, combining Lemma \ref{lsg1-lem} 
with a slight adaption of \reff{ntise} to the time 
dependent case we find that 
$C_1,C_2$ may be chosen independent of $n$. 
The higher regularity follows by a standard bootstrapping argument: 
for almost all $\tau\in (L^{-2},1)$ we have 
$(\uti^c_n,\uti^s_n)(\tau)\in B_{L^n}^3(3,2)$. 
Starting again at such a $\tau$ the 
required compatibility conditions to apply Lemma \ref{lsg1-lem} 
are automatically fulfilled. This yields \reff{hreg1}.\qed

%%%%%%%%%%%%%%%%%%%%%%%%%%%%%%%%%%%%%%%%%%%%%%%%%%%%%%%%%%%%%%%%%%%%%%%%%%%%

\subsection{Proof of Theorem \ref{thstabil-b} (Diffusive stability)}
\label{p1-sec}
Due to the loss of $L^{7/2}$ in Theorem \ref{loc-ex1} we need to 
improve \reff{ee1} to iterate \reff{r1}. Given a local solution 
$(\uti_n^c,\uti_n^s)$ with the higher regularity \reff{hreg1} this 
will be achieved by using the variation of constant formula and a 
suitable splitting of $\uti_n^c$. 

For $\uti^c\in \Hhat^3(2)$ we define $\Pi \uti^c= \pa_\kap \uti^c(0)$, 
which, by Sobolev embedding, gives a continuous map, 
i.e.\,$|\Pi \uti_c|\leq C\|\uti^c\|_{\Hhat^3(2)}$. 
Here, and also in \reff{lincontr1} below,  we need the smoothness in the 
Bloch wave number, which for $\uti^c$ we again identify with the Fourier 
wave number. To prove Theorem \ref{thstabil-b} we write 
\begin{gather} 
\begin{split}
\uti_n^c(\kap,1)&=\ri\psi_n g(\kap)+r_n^c(\kap), \qquad
\uti_n^s(\kap,\vt,1)=r_n^s(\vt,\kap),
\end{split}\label{ita}
\end{gather}
where $g(\kap)=\kap\er^{-\al \kap^2}$ and 
$r_n^c(0)=\pa_\kap r_n^c(0)=0$. This makes sense since $\uti_n(0,\tau)=0$ for 
all $n\in\N$ and all $\tau\in [1/L^2,1]$ if $\uti^c(0,1)=0$. 
Substituting \reff{ita} into \reff{splitn} yields 
\begin{subequations}
\label{it1}
\begin{align}
\psi_n-\psi_{n-1}&=\Pi I_n^c,\\
r_n^c&=\er^{(1-L^{-2})\lamti_{g,n}}L\CR_{1/L}r_{n-1}+I_n^c
+\res_n, \\
r_n^s&=\er^{(1-L^{-2})\Lamti_g}L\CR_{1/L}r_{n-1}+I_{n,0}^s+I_{n,1}^s+I_{n,2}^s,
\end{align}
\end{subequations}
where, using the notation $L^{3n}\CNti_n^s=\CNti_{n,0}^s+\CNti_{n,1}^s
+\CNti_{n,2}^s$ from Lemma \ref{non-lem}, 
\begin{align*}
I_n^c&=L^{3n}\int_{1/L^2}^1 \er^{(1-\tau)\lamti_{g,n}} 
\CNti_n^c(\uti_n^c,\uti_n^s)(\tau)\dd \tau,\quad 
I_{n,j}^s=\int_{1/L^2}^1 \er^{(1-\tau)\Lamti_{g,n}} 
\CNti_{n,j}^s(\uti_n^c,\uti_n^s)(\tau)\dd \tau, 
\end{align*}
and where the residual in (\ref{it1}b) is defined by 
$\res_n=\ri\psi_{n-1}\er^{(1-L^{-2})\lamti_{g,n}}L\CR_{1/L}g-\ri\psi_n g$. 
We also define 
$$
\rho_{n,c}=\|r_n^c\|_{B_{L^n}^3(2,2)}\quad\text{and}\quad 
\rho_{n,s}=\|r_n^s\|_{B_{L^n}^3(2,2)}, 
$$
which gives, cf.~\reff{rndef}, 
$\rho_n=\|\uti_n^c\|_{B_{L^n}^3(2,2)}+\|\uti_n^s\|_{B_{L^n}^3(2,2)}
\leq C|\psi_n|+\rho_{n,c}+\rho_{n,s}.$

Now assume that $\rho_{n-1}\leq L^{-7/2}$. Then from 
\reff{csge},\reff{ntice} we immediately obtain 
\begin{gather}
|\psi_n-\psi_{n-1}|\leq CL^{-n}(C_2 L^{7/2}\rho_{n-1})^2 
\label{psine1}
\end{gather} 
with $C_2$ from \reff{ee1}. Moreover,  
$\|\er^{(1-L^{-2})\lamti_{g,n}}L\CR_{1/L}g-g\|_{B^3_{L^n}(2,2)}
\leq CL^{-2n}$ and hence $\|\res_n\|_{B_{L^n}^3(2,2)}\leq CL^{-2n}|\psi_{n-1}|$. 
Next, we have 
\begin{gather}
\|\er^{(1-L^{-2})\lamti_{g,n}}L\CR_{1/L}
r_{n-1}^c\|_{B^3_{L^n}(2,2)}\leq CL^{-1}\|r_{n-1}^c\|_{B^3_{L^n}(2,2)}.
\label{lincontr1}
\end{gather}
This follows from 
$r_{n-1}^c(\kap/L)=(\frac{\kap}{L})^2\pa_\kap^2r_{n-1}^c(\tilde{\kap})$ for 
some $\tilde{\kap}$ between $0$ and $\kap$. Here again we need the 
smoothness in $\kap$. Combining the above estimates we arrive at 
\begin{gather}
\rho_{n,c}\leq CL^{-1}\rho_{n-1,c}+CL^{-n}(C_2L^{7/2}\rho_{n-1})^2+CL^{-2n}
|\psi_{n-1}|.  
\label{rhonce2}
\end{gather}

To estimate $\rho_{n,s}$ first note that
\huga{\label{mehr}
\|\er^{(1-L^{-2})\Lamti_g}L\CR_{1/L}r_{n-1}^s\|_{B^3_{L^n}(2,2)}
{\leq}C\er^{-\ga_0L^{2n}}L^{7/2}\rho_{n-1,s}{\leq}L^{-1}\rho_{n-1,s}
}
for $L$ sufficiently large. Next, 
$I_{n,0}^s$ and $I_{n,1}^s$ can be estimated using 
\reff{ssge} and \reff{ntise} to 
\begin{gather}
\|I_{n,0}^s\|_{B^3_{L^n}(2,2)}+\|I_{n,1}^s\|_{B^3_{L^n}(2,2)}
\leq CL^{-n}(L^{7/2}\rho_{n-1})^2. 
\label{rhonse1}
\end{gather}
However, for the quasi-linear part $I_{n,2}^s$ we have to use the 
higher regularity \reff{hreg1} and split 
$I_{n,2}^s=\int_{1/L^2}^{1/2}\cdots\dd\tau+\int_{1/2}^{1}\cdots\dd\tau$ 
to obtain 
\begin{align}
\|I_{n,2}^s\|_{B^3_{L^n}(2,2)}\leq& C(C_2L^{7/2}\rho_{n-1})^2\int_{1/L^2}^{1/2}
(1-\tau)^{-1}\er^{-\ga_0L^{2n}(1-\tau)}\dd\tau %\notag\\&
+C(C_3L^{7/2}\rho_{n-1})^2\int_{1/2}^{1}
\er^{-\ga_0L^{2n}(1-\tau)}\dd\tau\notag\\
\leq&C(C_2^2+C_3^2)L^{-n}(L^{7/2}\rho_{n-1})^2. 
\label{rhonse2}
\end{align}
This yields 
\huga{\label{rhonse2b}
\rho_{n,s}\le CL^{-1}\rho_{n-1}+CL^{-n}(L^{7/2}\rho_{n-1})^2
}
Now, let $L\geq L_0$ with $L_0$ so large that 
$CL^{-1}\leq L^{-(1-b)}$ and let 
$\rho_0=\|(\uti^c,\uti^s)\|_{B^3(2,2)}\leq L^{-4}$. 
Then, combining \reff{psine1}, \reff{rhonce2}, \reff{mehr} and 
\reff{rhonse2b}, iteration shows that there exists a $\psilim\in\R$ such 
that 
\begin{gather}
|\psilim-\psi_n|+\rho_{n,c}+\rho_{n,s}\leq L^{-n(1-b)}\quad\text{as}\quad 
n\ra\infty, 
\label{dcon}
\end{gather}
where the correction $L^{nb}$ takes care of the powers $C^n$ arising 
in the iteration. 
This discrete convergence implies Theorem \ref{thstabil-b} 
using $t=L^{2n}\tau$ and the local existence Theorem \ref{loc-ex1}. \qed 

%%%%%%%%%%%%%%%%%%%%%%%%%%%%%%%%%%%%%%%%%%%%%%%%%%%%%%%%%%%%%%%%%%%%%%%%%%%%

\subsection{Proof of Theorem \ref{thmix1-b}(i) (Diffusive mixing, 
Gaussian case)}\label{p2-sec}
The main difference compared to the proof of Theorem \ref{thstabil-b} 
are different scalings for $\uti_n^c,\uti_n^s$,  which are now based 
on \reff{rmap1} instead of \reff{rmap2}. Thus, we introduce 
\begin{gather} 
\uti_n^c(\kap,\tau)=\CR_{L^{-n}}\uti^c(\kap,L^{2n}\tau),\quad 
\uti_n^s(\vt,\kap,\tau)=\CR_{L^{-n}}\uti^s(\vt,\kap,L^{2n}\tau). 
\label{scal2}
\end{gather}
We want to show that 
\begin{gather}
\uti_n^c(\kap,1)\ra \phi_d\er^{-\al \kap^2},\quad 
\uti_n^s(\kap,1)\ra 0\quad\text{as}\quad n\ra\infty. 
\label{asy2}
\end{gather}
We obtain 
\begin{subequations}\label{splitn2}
\begin{align}
\pa_\tau \uti_n^c(\kap,\tau)-\lamti_{g,n}(\kap)\uti_n^c(\kap,\tau)
&=L^{2n}\CNti_n^c(\uti_n^c,\uti_n^s)( \kap,\tau),\label{unce2}\\
\pa_\tau \uti_n^s(\vt,\kap,\tau)-
\Lamti_{g,n}\uti_n^s(\vt,\kap,\tau)
&=L^{2n}\CNti_n^s(\uti_n^c,\uti_n^s)(\vt,\kap,\tau),\label{unse2}
\end{align}
\end{subequations}
where $\lamti_{g,n}(\kap)=L^{2n}\lamti_g(\kap/L^{2n})$ 
and $\Lamti_{g,n}=L^{2n}\CR_{L^{-n}}\Lamti_{g}\CR_{L^n}$ as before, and 
\begin{align*}
&\CNti_n^c(\uti_n^c,\uti_n^s)=
\etati(\kap/L^{n})\CR_{L^{-n}}\CNti(\CR_{L^n}\uti_n^c,\CR_{L^n}\uti_n^s), \\
&\CNti_n^s(\uti_n^c,\uti_n^s)=
\CR_{L^{-n}}\pmfsti\CNti(\CR_{L^n}\uti_n^c,\CR_{L^n}\uti_n^s).  
\end{align*}
The renormalization process reads 
\begin{gather} 
\text{solve \reff{splitn2} for } \tau\in I:=[L^{-2},1]\text{ with initial data }
\left.\bpm \uti_n^c\\\uti_n^s\epm\right|_{\tau=L^{-2}}=\CR_{1/L}
\left.\bpm\uti_{n-1}^c
\\ \uti_{n-1}^s\epm\right|_{\tau=1}.
\label{r3}
\end{gather}

The local existence for \reff{splitn2} works exactly as for \reff{splitn}. 
In contrast to Lemma \ref{loc-ex1} 
with $(\uti_n^c,\uti_n^s)(1/L^2)\in [B^3_{L^n}(2,2)]^2$ 
it suffices here to take $(\uti_n^c,\uti_n^s)(1/L^2)\in [B^2_{L^n}(2,2)]^2$ 
in order to extract the asymptotics \reff{asy2}, cf.\,\reff{lincontr2}. 
Thus, we set 
\begin{gather}
\rho_n:=\|(\uti_{n}^c,\uti_{n}^s)|_{\tau=1}\|_{[B_{L^n}^2(2,2)]^2},  
\label{rndef2}
\end{gather}
and for $\rho_{n-1}\leq C_1L^{-5/2}$ we obtain a local solution 
$(\uti_n^c,\uti_n^s)\in [\Kti_{L^n}^2(I,3,2)]^2$ to \reff{splitn2} with 
\begin{gather}
\|(\uti_n^c,\uti_n^s)\|_{[\Kti_{Ln}^2(I,3,2)]^2}\leq C_2L^{5/2}\rho_{n-1}, 
\label{ee1b}
\end{gather}
and for each $\tau_1>L^{-2}$ and $m_2\in\N$ there exists a $C_3$ such 
that we the higher regularity 
\begin{gather}
\|(\uti_n^c,\uti_n^s)\|_{[\Kti_{L^n}^2((\tau_1,1),m_2,2)]^2}
\leq C_3L^{5/2}\rho_{n-1}.
\label{hreg1b}
\end{gather}
The estimates for the linear semigroups work as before, i.e., here 
\begin{align}
\|\er^{(\tau-\tau')\lamti_{g,n}}\uti_n^c\|_{B_{L^n}^2(2,2)}
\leq &C\max\{1,(\tau-\tau')^{-b/2}\}\|\uti_n^c\|_{B_{L^n}^2(2,2-b)},
\label{csge2}\\
\|\er^{(\tau-\tau')\Lamti_{g,n}}\uti_n^s\|_{B_{L^n}^2(2,2)}
\leq &C\max\{1,(\tau-\tau')^{-m_2/2}\}\er^{-\ga_0L^{2n}(\tau-\tau')}
\|\uti_n^s\|_{B_{L^n}^2(2-m_2,2)}.
\label{ssge2}
\end{align}
The nonlinearities are now estimated as follows. 
\begin{Lemma}\label{non-lem2} 
Let $\|(\uti_n^c,\uti_n^s)\|_{[B^2_{L^n}(2,2)]^2}\leq R_n\leq 1$. 
There exists a $C>0$ such that 
\begin{align}
L^{2n}\|\CNti_n^c(\uti_n^c,\uti_n^s)\|_{B_{L^n}^2(2,1)}&\leq CL^{-n}R_n^2. 
\label{ntice2}
\end{align}
The term
$L^{2n}\CNti_n^s$ can be split according to the number of $\vt$ derivatives 
as $L^{2n}\CNti_n^s=\CNti_{n,0}^s+\CNti_{n,1}^s+\CNti_{n,2}^s$ such that 
\begin{align}
\|\CNti_{n,i}^s\|_{B_{L^n}^2(2-i,2)}&\leq CR_n^2.
\label{ntise2}
\end{align}
\end{Lemma}
{\bf Proof.} Apart from the different power counting, the proof works 
like the one of Lemma \ref{non-lem}, with the crucial difference that 
now the term $s_1$ from \reff{s1def} vanishes since $\beta=0$ by assumption. 
\qed 

Similar to \reff{ita} we now set 
\begin{gather} 
\begin{split}
\uti_n^c(\kap,1)&=\phi_d g(\kap)+r_n^c(\kap), \qquad
\uti_n^s(\kap,\vt,1)=r_n^s(\kap,\vt),
\end{split}\label{ita2}
\end{gather}
where $g(\kap)=\er^{-\al \kap^2}$ and $r_n^c(0)=0$. Here no 
variables $\phi_n$ are necessary since, due to the conservation 
of total phase shift, i.e., $\pa_t\uti^c(0,t)=0$ for 
all $t$. We obtain 
\begin{subequations}
\label{it2}
\begin{align}
r_n^c&=\er^{(1-L^{-2})\lamti_{g,n}}\CR_{1/L}r_{n-1}^c+I_n^c+\res_n, \\
r_n^s&=\er^{(1-L^{-2})\Lamti_g}\CR_{1/L}r_{n-1}^s+I_{n,0}^s+I_{n,1}^s+I_{n,1}^s,
\end{align}
\end{subequations}
where $\res_n=\phi_d(\er^{(1-L^{-2})\lamti_{g,n}}\CR_{1/L}g-g)$ and 
\begin{align*}
I_n^c&=L^{2n}\int_{1/L^2}^1 \er^{(1-\tau)\lamti_{g,n}} 
\CNti_n^c(\uti_n^c,\uti_n^s)(\tau)\dd \tau,\quad 
I_{n,j}^s=\int_{1/L^2}^1 \er^{(1-\tau)\Lamti_{g,n}} 
\CNti_{n,j}^s(\uti_n^c,\uti_n^s)(\tau)\dd \tau, 
\end{align*}
with $L^{2n}\CNti_n^s=\CNti_{n,0}^s+\CNti_{n,1}^s
+\CNti_{n,2}^s$ from Lemma \ref{non-lem2}. 
Clearly, $\|\res_n\|_{B_{L^n}^2(2,2)}\leq C|\phi_d|L^{-2n}$, and 
\begin{gather}
\|\er^{(1-L^{-2})\lamti_{g,n}}\CR_{1/L}
r_{n-1}^c\|_{B^2_{L^n}(2,2)}\leq CL^{-1}\|r_{n-1}^c\|_{B^2_{L^n}(2,2)} 
\label{lincontr2}
\end{gather}
which follows by writing 
$r_{n-1}^c(\kap/L)=(\frac{\kap}{L})\pa_\kap r_{n-1}^c(\tilde{\kap})$ for 
some $\tilde{\kap}$ between $0$ and $\kap$. Combining this 
with \reff{csge2} and \reff{ntice2} thus yields 
\begin{gather}
\rho_{n,c}\leq CL^{-1}\rho_{n-1,c}+CL^{-n}(C_2L^{5/2}\rho_{n-1})^2 
+C|\phi_d|L^{-2n}, 
\label{rhonce3}
\end{gather}
and similarly 
\begin{gather}
\rho_{n,s}\leq L^{-1}\rho_{n-1,s}+C(C_2^2+C_3^2)L^{-n}(L^{5/2}\rho_{n-1})^2. 
\label{rhonse3}
\end{gather}
The proof of Theorem \ref{thmix1-b}(i) 
now follows by iteration. At this point the assumption $|\phi_d|\le \eps$ 
is crucial. This can be seen by computing $ \rho_{n,c} $ in powers 
of $ L $ for $ n = 0,1,2, \ldots $ starting 
with $ \rho_{0,c}=0$. Hence, we need $ |\phi_d| \leq  L^{-d} $ 
with $ d $ sufficiently large.
\qed

%%%%%%%%%%%%%%%%%%%%%%%%%%%%%%%%%%%%%%%%%%%%%%%%%%%%%%%%%%%%%%%%%%%%%%%%%%%%

\subsection{Proof of Theorem \ref{thmix1-b}(ii) (Diffusive Mixing, Burgers' case)}\label{p3-sec}

Essentially, Theorem \ref{thmix1-b}(ii) is again based on the 
scaling \reff{scal2}. However, the crucial difference to the case 
$\beta=0$ is that now the analog of \reff{ntice2} no longer holds. 
Therefore, we need to scale $\uti^s$ differently, i.e., we blow up 
$\uti_n^s$ in order to avoid problems with the quadratic terms 
involving $\uti^s$ in the critical part, in which the term 
$\ri\ell\beta(\uti^c*\uti^c)$ 
will give the Burgers dynamics for $\uti^c$. Thus, 
for small $p>0$ we introduce 
\begin{gather} 
\uti_n^c(\kap,\tau)=\CR_{L^{-n}}\uti^c(\kap,L^{2n}\tau),\quad 
\uti_n^s(\vt,\kap,\tau)=L^{n(1-p)}\CR_{L^{-n}}\uti^s(\vt,\kap,L^{2n}\tau), 
\label{scal3}
\end{gather}
to obtain 
\begin{subequations}\label{splitn3}
\begin{align}
\pa_\tau \uti_n^c(\kap,\tau)-\lamti_{g,n}(\kap)\uti_n^c(\kap,\tau)
&=L^{2n}\CNti_n^c(\uti_n^c,\uti_n^s)( \kap,\tau),\label{unce3}\\
\pa_\tau \uti_n^s(\vt,\kap,\tau)-
\Lamti_{g,n}\uti_n^s(\vt,\kap,\tau)
&=L^{n(3-p)}\CNti_n^s(\uti_n^c,\uti_n^s)\vt,\kap,\tau),\label{unse3}
\end{align}
\end{subequations}
where again $\lamti_{g,n}(\kap)=L^{2n}\lamti_g(\kap/L^{n})$ 
and $\Lamti_{g,n}=L^{2n}\CR_{L^{-n}}\Lamti_{g}\CR_{L^n}$, but now 
\begin{align*}
&\CNti_n^c(\uti_n^c,\uti_n^s)=
\etati(\kap/L^{n})\CR_{L^{-n}}\CNti(\CR_{L^n}\uti_n^c,
L^{-n(1-p)}\CR_{L^n}\uti_n^s), \\
&\CNti_n^s(\uti_n^c,\uti_n^s)=
\CR_{L^{-n}}\pmfsti\CNti(\CR_{L^n}\uti_n^c,L^{-n(1-p)}\CR_{L^n}\uti_n^s). 
\end{align*}
Accordingly, the renormalization process reads 
\begin{gather} 
\begin{split}
\text{solve \reff{splitn3} on } \tau\in I:=[L^{-2},1]\text{ with initial data } 
\left.\bpm \uti_n^c\\\uti_n^s\epm\right|_{\tau=L^{-2}}=\CR_{1/L}
\left.\bpm\uti_{n-1}^c
\\ L^{1-p}\,\uti_{n-1}^s\epm\right|_{\tau=1}. 
\end{split}
\label{r2}
\end{gather}
The estimates for the linear semigroups are again \reff{csge2} and \reff{ssge2}, 
while the nonlinear terms are estimated as follows.
\begin{Lemma}\label{non-lem3} 
Let $\|(\uti_n^c,\uti_n^s)\|_{[B^2_{L^n}(2,2)]^2}\leq R_n\leq 1$. 
There exists a $C>0$ such that 
$L^{2n}\CNti_n^c(\uti_n^c,\uti_n^s)=s_1+s_2+s_3+s_4$ with 
$s_1=\ri \beta\kap(\uti_n^c*\uti_n^c)$ and 
\begin{gather}
\begin{split}
\|s_2\|_{B^2_{L^n}(2,1)}&\leq CL^{-n}\|\uti_n^c\|_{B^2_{L^n}(2,2)}^2,\\
\|s_3\|_{B^2_{L^n}(2,1)}&\leq CL^{-n(1-p)}\|\uti_n^c\|_{B^2_{L^n}(2,2)}
\|\uti_n^s\|_{B^2_{L^n}(2,2)}, \\
\|s_4\|_{B^2_{L^n}(2,1)}&\leq CL^{-2n(1-p)}R_n^2. 
\end{split}\label{ntice3}
\end{gather}
The term
$L^{n(3-p)}\CNti_n^s$ can be split according to the number of $\vt$ derivatives 
as $L^{n(3-p)}\CNti_n^s=\CNti_{n,0}^s+\CNti_{n,1}^s+\CNti_{n,2}^s$ such that 
\begin{align}
\|\CNti_{n,i}^s\|_{B_{L^n}^2(2-i,2)}&\leq C(L^{n(1-p)}
\|\uti_n^c\|_{B_{L^n}^2(2,2)}^2+R_n^2).
\label{ntise3}
\end{align}
\end{Lemma}
{\bf Proof.} 
The term $s_2$ contains the quadratic terms in $\uti_n^c$ except 
for $\ri \beta\kap(\uti_n^c*\uti_n^c)$, i.e., $s_2$ is of the form 
$s_2(\kap,t)=h(\kap/L^n)(\uti_n^c*\uti_n^c)$ with $h(\kap)=\CO(\kap^2)$. 
$s_3$ contains the quadratic interaction of $\uti_n^c$ and $\uti_n^s$, and 
$s_4$ contains the remaining terms. Then \reff{ntice3} follows from 
% \reff{bcone1} and by applying Lemma \ref{dxlem} and taking into account 
the finite support of $\CNti_n^c$ in Fourier space. It is in $s_3,s_4$ that 
the blowup scaling 
$\uti_n^s(\vt,\kap,\tau)=L^{n(1-p)}\CR_{L^{-n}}\uti^s(\vt,\kap,L^{2n}\tau)$ 
is useful. Equation \reff{ntise3} again 
follows by straightforward power counting. 
\qed 

The terms involving only $\uti_n^c$ in \reff{ntise3} blow up 
as $n\ra\infty$. However, combining \reff{ntise3} with the exponential 
damping in the stable part, cf.\,\reff{ssge2}, we still get local 
existence for \reff{splitn3} 
with constants independent of $n$. We let 
\begin{gather}
\rho_n:=\|(\uti_{n}^c,\uti_{n}^s)|_{\tau=1}\|_{[B_{L^n}^2(2,2)]^2},  
\label{rndef3}
\end{gather}
and for $\rho_{n-1}\leq C_1L^{-5/2}$ obtain a local solution 
$(\uti_n^c,\uti_n^s)\in [\Kti_n^2(I,3,2)]^2$ to \reff{splitn3} with 
$\|(\uti_n^c,\uti_n^s)\|_{[\Kti_{L^n}^2(I,3,2)]^2}\leq C_2L^{5/2}\rho_{n-1}$, 
as in \reff{ee1b}, 
% \begin{gather}
% \|(\uti_n^c,\uti_n^s)\|_{[\Kti_n^2(I,3,2)]^2}\leq C_2L^{5/2}\rho_{n-1}, 
% \label{ee1c}
% \end{gather}
which moreover enjoys the higher regularity 
$\|(\uti_n^c,\uti_n^s)\|_{[\Kti_{L^n}^2((\tau_1,1),m_2,2)]^2}
 \leq C_3L^{5/2}\rho_{n-1}$, cf.~\reff{hreg1b}. 
% \begin{gather}
% \|(\uti_n^c,\uti_n^s)\|_{[\Kti_n^2((\tau_1,1),m_2,2)]^2}
% \leq C_3L^{5/2}\rho_{n-1}.
% \label{hreg1c}
% \end{gather}

Similar to \reff{ita} and \reff{ita2} we now separate from $\uti_n^c$ 
the lowest order asymptotics, now obtained from the Burgers equation. 
However, due to the contribution 
$s_1=\ri\beta\kappa(\uti_n^c*\uti_n^c)$ of the nonlinearity 
to the asymptotics here we work 
out an intermediate step and split $\uti_n^c$ 
in a $\tau$ dependent way. 
In detail, let 
\begin{gather} 
\uti_n^c(\kap,\tau)=\uti_{n,*}^c(\kap,\tau)+\alti_n(\kap,\tau)
\label{uncsplit5}
\end{gather}
where 
$$\uti_{n,*}^c(\kap,\tau)=\chi(\kap/L^n)\uti_{*}^c(\kap,\tau)
\quad\text{with}\quad
\uti_*^c(\ell,t)=
\CF\left(\frac{\sqrt{\al}}\beta\frac{z\er^{-\vt^2/(k^2\al t)}}
{1+z\erf(\vt/\sqrt{k\al t})}\right)(\ell), 
$$
cf.~\reff{n27c}, 
and $\chi$ from \reff{chidef}. Consequently 
$\alti_n(0,\tau){=}0$ for all $n,\tau$ due to the conservation of 
total phase.  
Then 
\begin{subequations}\label{splitn4}
\begin{align}
\pa_\tau \alti_n^c(\kap,\tau)-\lamti_{g,n}(\kap)\alti_n^c(\kap,\tau)
&=L^{2n}(\CNti_n^c(\uti_n^c,\uti_n^s)( \kap,\tau)-
\CNti_n^c(\uti_{n,*}^c,0))+\res_n,\label{unce4}\\
\pa_\tau \uti_n^s(\vt,\kap,\tau)-
\Lamti_{g,n}\uti_n^s(\vt,\kap,\tau)
&=L^{n(3-p)}\CNti_n^s(\uti_n^c,\uti_n^s)(\vt,\kap,\tau),\label{unse4}
\end{align}
\end{subequations}
where 
$$
\res_n=-\pa_\tau\uti_{n,*}^c+\lamti_{g,n}(\kap)\uti_{n,*}^c
+L^{2n}\CNti_n^c(\uti_{n,*}^c,0). 
$$
\begin{Lemma}\label{res-lem}
There exists a $C>0$ such that 
$\sup_{L^{-2}\leq \tau\leq 1}
\|\res_n(\tau)\|_{B_{L^n}^2(2,2)}\leq CL^{-n}|\phi_d|$. 
\end{Lemma}
{\bf Proof.}\ By construction, i.e., since $\uti_{n,*}^c$ 
is an exact solution of the Burgers equation, 
$$
\res_n(\kap,\tau)=CL^{-n}\biggl(\CO(\kap^3)\uti_{n,*}^c+
\CO(\kap^2(\uti_{n,*}^c*\uti_{n,*}^c))\biggr)
$$
which can be estimated in $B_{L^n}^2(2,2)$ by $L^{-n}|\phi_d|$ since 
$\uti_{*}^c(\kap,\tau)$ is analytic and exponentially decaying.\qed

Now setting 
\begin{gather} 
\begin{split}
\uti_n^c(\kap,1)&=\uti_{n,*}^c(\kap,1)+r_n^c(\kap), \qquad
\uti_n^s(\kap,\vt,1)=r_n^s(\kap,\vt),
\end{split}\label{ita3}
\end{gather}
the remainder of the proof of Theorem \ref{thmix1-b}(ii) works as 
the proof of Theorem \ref{thmix1-b}(i) in \S\ref{p2-sec}. \qed

%%%%%%%%%%%%%%%%%%%%%%%%%%%%%%%%%%%%%%%%%%%%%%%%%%%%%%%%%%%%%%%%%%%%%%%%%%%%
%\newpage
%\bibliographystyle{alpha}\bibliography{/home/hannes/texb/hubib}

\begin{thebibliography}{DSSS06}

\bibitem[BK92]{bk92}
J.~Bricmont and A.~Kupiainen.
\newblock Renormalization group and the {G}inzburg--{L}andau equation.
\newblock {\em Comm. Math. Phys.}, 150:193--208, 1992.

\bibitem[BKL94]{bkl94}
J.~Bricmont, A.~Kupiainen, and G.~Lin.
\newblock Renormalization group and asymptotics of solutions of nonlinear
  parabolic equations.
\newblock {\em Comm. Pure Appl. Math.}, 6:893--922, 1994.

\bibitem[DSSS09]{DSSS05}
A.~Doelman, B.~Sandstede, A.~Scheel, and G.~Schneider.
\newblock The dynamics of modulated wave trains.
\newblock {\it {M}emoirs of the {A}{M}{S}}, 199(934), 2009.

\bibitem[ES00]{es00b}
J.-P. Eckmann and G.~Schneider.
\newblock Nonlinear stability of bifurcating front solutions for the
  {T}aylor--{C}ouette problem.
\newblock {\em Special Issue of ZAMM}, 80(11--12):745--753, 2000.

\bibitem[GM98]{gm98}
Th. Gallay and A.~Mielke.
\newblock Diffusive mixing of stable states in the {G}inzburg--{L}andau
  equation.
\newblock {\em Comm. Math. Phys.}, 199:71--97, 1998.

\bibitem[GSU04]{gsu04}
Th. Gallay, G.~Schneider, and H.~Uecker.
\newblock Stable transport of information near essentially unstable localized
  structures.
\newblock {\em DCDS--B}, 4(2):349--390, 2004. 

\bibitem[JNRZ11]{jnrz11}
M.~Johnson, P.~Noble, L.~M. Rodrigues, and K.~Zumbrun.
\newblock Nonlocalized modulation of periodic reaction diffusion waves:
  Nonlinear stability.
\newblock Preprint, 2011.

\bibitem[JZ11]{jz10}
M.~Johnson and K.~Zumbrun.
\newblock Nonlinear stability of spatially-periodic traveling-wave solutions of
  systems of reaction diffusion equations.
\newblock {\it Annales de l'Institut Henri Poincare - Analyse Non Lineaire}, to
  appear, 2011.

\bibitem[LM68]{lm68}
J.L. Lions and E.~Magenes.
\newblock {\em Probl\`emes aux limites non homog\`enes}.
\newblock Dunod, Paris, 1968.

\bibitem[MSU01]{msu01}
A.~Mielke, G.~Schneider, and H.~Uecker.
\newblock Stability and diffusive dynamics on extended domains.
\newblock In {\em Ergodic theory, analysis, and efficient simulation of
  dynamical systems}, pages 563--583. Springer, Berlin, 2001.

\bibitem[RS78]{RS72}
M.~Reed and B.~Simon.
\newblock {\em Methods of Modern Mathematical Physics IV}.
\newblock Academic Press, 1978.

\bibitem[Sca99]{scar99}
B.~Scarpellini.
\newblock {\em Stability, Instability, and Direct Integrals}.
\newblock Chapman \& Hall, 1999.

\bibitem[Sch94]{Schn94ZAMP}
G.~Schneider.
\newblock Error estimates for the {G}inzburg--{L}andau approximation.
\newblock {\em ZAMP}, 45:433--457, 1994.

\bibitem[Sch96]{gs96a}
G.~Schneider.
\newblock Diffusive stability of spatial periodic solutions of the
  {S}wift--{H}ohen\-berg equation.
\newblock {\em Comm. Math. Phys.}, 178:679--702, 1996.

\bibitem[Sch98a]{gstohoku}
G.~Schneider.
\newblock Nonlinear diffusive stability of spatially periodic
  solutions---abstract theorem and higher space dimensions.
\newblock In {\em Proceedings of the International Conference on Asymptotics in
  Nonlinear Diffusive Systems (Sendai, 1997)}, volume~8 of {\em Tohoku Math.
  Publ.}, pages 159--167, 1998.

\bibitem[Sch98b]{gs98}
G.~Schneider.
\newblock Nonlinear stability of {T}aylor--vortices in infinite cylinders.
\newblock {\em Arch. Rat. Mech. Anal.}, 144(2):121--200, 1998.

\bibitem[Uec04]{ue04ibl}
H.~Uecker.
\newblock Self--similar decay of localized perturbations in the {I}ntegral
  {B}oundary {L}ayer equation.
\newblock {\em JDE}, 207(2):407--422, 2004.

\bibitem[Uec07]{ue04nsb}
H.~Uecker.
\newblock Self-similar decay of spatially localized perturbations of the
  {N}usselt solution for the inclined film problem.
\newblock {\em Arch.\,Rat.\,Mech.\,Anal.}, 184(3):401--447, 2007.

\end{thebibliography}
\begin{small}
\end{small}

\end{document}